\documentclass[11pt,a4paper]{article}

\usepackage{header}
\usepackage{abbrev}

\begin{document}

\title{Implicit low-rank Riemannian schemes for the time integration of stiff partial differential equations}

\author{Marco Sutti\thanks{Mathematics Division, National Center for Theoretical Sciences, Taipei, Taiwan (\email{msutti@ncts.tw}).}\hspace{2mm}\orcidlink{0000-0002-8410-1372}
~and Bart Vandereycken\thanks{Section of Mathematics, University of Geneva, Geneva, Switzerland (\email{bart.vandereycken@unige.ch}).}}

\date{\today}

\maketitle

\begin{abstract}
We propose two implicit numerical schemes for the low-rank time integration of stiff nonlinear partial differential equations. Our approach uses the preconditioned Riemannian trust-region method of Absil, Baker, and Gallivan, 2007. We demonstrate the efficiency of our method for solving the Allen--Cahn and the Fisher--KPP equation on the manifold of fixed-rank matrices. Furthermore, our approach allows us to avoid the restriction on the time step typical of methods that use the fixed-point iteration to solve the inner nonlinear equations.
Finally, we demonstrate the efficiency of the preconditioner on the same variational problems presented in Sutti and Vandereycken, 2021.

\bigskip
\textbf{Key words.} Implicit methods, numerical time integration, Riemannian optimization, stiff PDEs, manifold of fixed-rank matrices, variational problems, preconditioning, trust-region method, Allen--Cahn equation, Fisher--KPP equation

\medskip
\textbf{AMS subject classifications.} 65F08, 65F55, 65L04, 65F45, 65N22, 65K10, 58C05

\end{abstract}

\section{Introduction}

The topic of this paper is the efficient solution of large-scale variational problems arising from the discretization of partial differential equations (PDEs), both time-independent and time-dependent. 
In the first part of the paper, we use the preconditioned Riemannian trust-region method of Absil, Baker, and Gallivan~\protect{\cite{ABG:2007}} to solve the nonlinear equation derived from an implicit scheme for numerical time integration. All the calculations are performed on a low-rank manifold, which allows us to approximate the solution with significantly fewer degrees of freedom.
In the second part of the paper, we solve variational problems derived from the discretization of elliptic PDEs. These are large-scale finite-dimensional optimization problems arising from the discretization of infinite-dimensional problems. Variational problems of this type have been considered as benchmarks in several nonlinear multilevel algorithms \protect{\cite{Henson:2003,Gratton:2008,Wen:2009,Sutti_V:2021}}.

A common way to speed up numerical computations is by approximating large matrices using low-rank methods. This is particularly useful for high-dimensional problems, which can be solved using low-rank matrix and tensor methods. Tensors are simply the higher dimensional version of two-dimensional matrices, as explained in Hackbusch's work~\cite{Hackbusch:2012}. The earliest examples are low-rank solvers for the Lyapunov equation, $AX+XA\tr = C$, and other matrix equations, see, e.g., \protect{\cite{KressnerTobler:2011,Simoncini:2016,Kurschner:2016}}. The low-rank approximation properties for these problems are also reasonably well-understood from a theoretical point of view. Indeed, \cite[Remark~1]{Grasedyck:2004} showed that the solution $ X $ to a Sylvester equation $ AX - XB + C = 0 $ could be approximated up to a relative accuracy of $ \varepsilon $ using a rank $ r = \cO(\log_{2}(\mathrm{cond}_{2}(A)) \log_{2}(1/\varepsilon) ) $.
Typically, to obtain a low-rank approximation of the unknown solution $X$, an iterative method that directly constructs the low-rank approximation is used. This work falls under the category of techniques that achieve this through \emph{Riemannian optimization}~\cite{Edelman:1998,AMS:2008,boumal_2023}. To ensure critical points have a low-rank representation, the optimization problem (which may be reformulated from the original) is limited to the fixed-rank matrices manifold $\cMr$. Some early references on this manifold include \protect{\cite{Helmke:1994,Koch:2007,Vandereycken:2013}}. Retraction-based optimization on $ \cMr $ was studied in \protect{\cite{Shalit:2010,Shalit:2012}}. Optimization on $ \cMr $ has gained a lot of momentum during the last decade, and examples of such methods are~\cite{Mishra:2011b,Vandereycken:2013,Steinlechner:2016} for matrix and tensor completion, \cite{Shalit:2012} for metric learning, \cite{VandereyckenV_2010,Mishra_V_2014,KSV:2016} for matrix and tensor equations, and \cite{Rakhuba:2018ab,RAKHUBA2019718} for eigenvalue problems. These optimization problems are ill-conditioned in discretized PDEs, making simple first-order methods like gradient descent excessively slow.

\subsection{Riemannian preconditioning}
This paper is closely related to the concept of preconditioning on Riemannian manifolds, which is similar to preconditioning in the unconstrained case (see, e.g., \protect{\cite{NW:2006}}). Several authors have tackled preconditioning in the Riemannian optimization framework; the following overview is not meant to be exhaustive. In \cite{VandereyckenV_2010,KSV:2016,Rakhuba:2018ab}, for example, the gradient is preconditioned with the inverse of the local Hessian. Solving these Hessian equations is done by a preconditioned iterative scheme, mimicking the class of quasi or truncated Newton methods. We also refer to~\cite{UschmajewV:2020} for a recent overview of geometric methods for obtaining low-rank approximations. The work most closely related to the present paper is~\protect{\cite{VandereyckenV_2010}}, which proposed a preconditioner for the manifold of symmetric positive semidefinite matrices of fixed rank.
\protect{\cite{Boumal_Absil:2015}} developed a preconditioner for Riemannian optimization on the Grassmann manifold. \protect{\cite{Mishra:2016}} investigated the connection between quadratic programming and Riemannian gradient optimization, particularly on quotient manifolds. Their method proved efficient, especially in quadratic optimization with orthogonality and rank constraints. Related to this preconditioned metric approach are those of~\protect{\cite{Ngo_Saad:2012,Mishra_V_2014}}, and more recently~\protect{\cite{cai2022tensor}}, who extend the preconditioned metric from the matrix case to the tensor case using the tensor train (TT) format for the tensor completion problem. On tensor manifolds, \protect{\cite{KSV:2016}} developed a preconditioned version for Riemannian gradient descent and the Richardson method, using the Tucker and TT formats.

\subsection{Trust-region methods}
This paper uses the Riemannian trust-region method of Absil, Baker, and Gallivan~\protect{\cite{ABG:2007}}. This method embeds an inner truncated conjugate gradient (tCG) method to solve the so-called trust-region minimization subproblem. The tCG solver naturally lends itself to preconditioning, and the preconditioner is typically a symmetric positive definite operator that approximates the inverse of the Hessian matrix. Ideally, it has to be cheap to compute.
Preconditioning with the projected Euclidean Hessian was done for symmetric positive semidefinite matrices with fixed rank~\cite{VandereyckenV_2010}. In contrast, we develop it here for any, i.e., typically non-symmetric, fixed-rank matrix.

We follow the steps outlined in~\protect{\cite{VandereyckenV_2010}} to find the preconditioner, namely: find the Euclidean Hessian, find the Riemannian Hessian operator, vectorize it to get the Hessian matrix, a linear and symmetric matrix; the inverse of the Hessian matrix should make a good candidate for a preconditioner; apply the preconditioner.

\subsection{Low-rank approximations for time-dependent PDEs}
Various approaches have been employed to address the low-rank approximation of time-dependent partial differential equations (PDEs). One such method is the dynamical low-rank approximation (DLRA)~\protect{\cite{Koch:2007,Lubich:2014}}, which optimally evolves a system’s low-rank approximation for common time-dependent PDEs. For example, suppose we are given a discretized dynamical system as a first-order differential equation (ODE). The DLRA idea is to replace the $ \dt{W} $ in the ODE with the tangent vector in $ \mathrm{T}_{W}\cMr $ that is closer to the right-hand side $ G(W) $. Recent developments of the DLRA include, but are not limited to, \protect{\cite{KieriLubichWalach:2016,KieriVandereycken:2019,Musharbash:2020,CerutiLubich:2022,CerutiKuschLubich:2022}}, and \protect{\cite{Billaud-Friess:2022}}.

Another approach is the dynamically orthogonal Runge--Kutta of~\protect{\cite{LermusiauxRobinson:1999,Lermusiaux:2001}}, and its more recent developments \protect{\cite{SapsisLermusiaux:2009,Ueckermann:2013,FepponLermusiaux:2018aa,FepponLermusiaux:2018bb,CharousLermusiaux:2021,charous2022stable}}.

The step-truncation methods of~\protect{\cite{RodgersDektorVenturi:2022,RodgersVenturi:2022}} form another class of methods for the low-rank approximation of time-dependent problems.
In~\protect{\cite{RodgersVenturi:2022}}, they study implicit rank-adaptive algorithms based on performing one time step with a conventional time-stepping scheme, followed by an implicit fixed-point iteration step involving a rank truncation operation onto a tensor or matrix manifold.

Here, we follow another route. While also employing an implicit time-stepping scheme for the time evolution, as in~\protect{\cite{RodgersVenturi:2022}}, instead of using a fixed-point iteration method for solving the nonlinear equations, we use a preconditioned Riemannian trust-region method (named PrecRTR) on the manifold of fixed-rank matrices. 
This results in a preconditioned dynamical low-rank approximation of the Allen--Cahn and Fisher--KPP equations.

Although implicit time integration methods are more expensive, they allow for a larger time step than their explicit counterparts. Moreover, the cost of solving the inner nonlinear equations remains moderately low thanks to our preconditioner.

Recently, \protect{\cite{Massei:2022}} also investigated the low-rank numerical integration of the Allen--Cahn equation, but their approach is very different from ours since they use hierarchical low-rank matrices.

\subsection{Contributions and outline}
The most significant contribution of this paper is an implicit numerical time integration scheme that can be used to solve stiff nonlinear time-dependent partial differential equations (PDEs). Our method internally employs a preconditioned Riemannian trust-region method on the manifold of fixed-rank matrices to solve the implicit equation derived by the time integration scheme.
We also include the development of a preconditioner for the Riemannian trust-region (RTR) subproblem on the manifold of fixed-rank matrices. This can be regarded as an extension of the preconditioner of~\protect{\cite{VandereyckenV_2010}} for the manifold of symmetric positive semidefinite matrices of fixed rank. We focus on applying the preconditioned RTR to the solution of two time-dependent, stiff nonlinear PDEs: the Allen--Cahn and Fisher--KPP equations. Additionally, we consider the two variational problems already studied in~\protect{\cite{Henson:2003,Gratton:2008,Wen:2009,Sutti_V:2021}}. The numerical experiments demonstrate the efficiency of the preconditioned algorithm in contrast to the non-preconditioned algorithm.

The remaining part of this paper is organized as follows. Section~\ref{sec:problem_setting} introduces the problem settings and the objective functions of the problems, which are the focus of this work. In section~\ref{sec:riem_optim}, we recall some preliminaries on the Riemannian optimization framework and the Riemannian trust-region method. Section~\ref{sec:geom_fixed_rank_manif} details the geometry of the manifold of fixed-rank matrices, tangent spaces, projectors, Riemannian gradient and Hessian, and the retraction.  In section~\ref{sec:RTR}, we recall more algorithmic details of the Riemannian trust-region method. Section~\ref{sec:ACE} and~\ref{sec:FKPP} present the core contribution of this paper: an implicit Riemannian low-rank scheme for the numerical integration of stiff nonlinear time-dependent PDEs, the Allen--Cahn equation and the Fisher--KPP equation. Other numerical experiments on the two variational problems from~\protect{\cite{Sutti_V:2021}} are presented and discussed in section~\ref{sec:numerical_experiments}. Finally, we wrap up our work in section~\ref{sec:conclusions}. More details about the derivation of the preconditioner for the Riemannian trust-region method on the manifold of fixed-rank matrices are given in Appendix~\ref{app:preconditioner}, and the discretization details for the Allen--Cahn and the Fisher--KPP equations are provided in Appendices~\ref{app:discretization_ACE} and~\ref{app:discretization_FKPP}, respectively.

\subsection{Notation}
The space of $ n \times r $ matrices is denoted by $ \R^{n \times r} $. By $ X_{\perp} \in \R^{n \times (n-r)}$ we denote an orthonormal matrix whose columns span the orthogonal complement of $ \mathrm{span}(X) $, i.e., $ X\tr_{\perp} X = 0 $ and $ X_{\perp}\tr X_{\perp} = I_{r} $.
In the formulas throughout the paper, we typically use the Roman capital script for operators and the italic capital script for matrices. For instance, $ \PTXM $ indicates a projection \emph{operator}, while $\PXmat $ is the corresponding projection \emph{matrix}.

The directional derivative of a function $ f $ at $ x $ in the direction of $ \xi $ is denoted by $ \D\! f(x)[\xi] $. With $ \| \cdot \|_{\F} $, we indicate the Frobenius norm of a matrix.

Even though we did not use a multilevel algorithm in this work, we want to maintain consistency with the notation used in~\protect{\cite{Sutti_V:2021}}. Consequently, we denote by $ \ell $ the discretization level. Hence the total number of grid points on a two-dimensional square domain is given by $ 2^{2\ell} $. This notation was adopted in~\protect{\cite{Sutti_V:2021}} due to the multilevel nature of the Riemannian multigrid line-search (RMGLS) algorithm. Moreover, here we omit the subscripts $ \cdot_{h} $ and $ \cdot_{H} $ because they were due to the multilevel nature of RMGLS. We use $ \Delta $ to denote the Laplacian operator, and the spatial discretization parameter is denoted by $ h_{x} $. The time step is represented by $ \deltatime $.

\section{The problem settings and cost functions} \label{sec:problem_setting}
In this section, we present the optimization problems studied in this paper. The first two problems are time-dependent, stiff PDEs: the Allen--Cahn and thes the Fisher--KPP equations.
The last two problems are the same considered in~\protect{\cite{Sutti_V:2021}}.

The Allen--Cahn equation in its simpler form reads (see Sect.~\ref{sec:ACE})
\[
   \frac{\partial w}{\partial t} = \varepsilon \Delta w + w  - w^{3},
\]
where $ w \equiv w(\bm{x}, t) $, $ \bm{x} \in \Omega = [-\pi, \pi]^{2} $, and $ t\geq 0 $. We solve it on a two-dimensional flat torus.
We reformulate it as a variational problem, which leads us to consider
\[
   \min_{w} \cF(w) \coloneqq \int_{\Omega} \frac{\varepsilon \deltatime}{2} \| \nabla w \|^{2} + \frac{(1-\deltatime)}{2} \, w^{2} + \frac{\deltatime}{4} \, w^{4} - \widetilde{w} \cdot w \dx \dy.
\]

The second problem considered is the Fisher--KPP equation with homogeneous Neumann boundary conditions, for which we construct the cost function (see Sect.~\ref{sec:FKPP})
\begin{small}
\[
  F(W) = \frac{1}{2} \trace\!\big( W\tr \! M_{\mathrm{m}}\tr \! M_{\mathrm{m}} W \big) - \trace\!\big((W^{(n-1)})\tr \! M_{\mathrm{p}}\tr \! M_{\mathrm{m}} W\big) + 2 \deltatime \trace\!\left(\left(\big( W^{(n)} \big)^{\circ 2} -  W^{(n)}  \right)\tr \! M_{\mathrm{m}} W R_{\omega}\right).
\]
\end{small}

Thirdly, we study the following variational problem, studied in~\protect{\cite{Henson:2003,Gratton:2008,Wen:2009}}, and called ``LYAP'' in~\protect{\cite[Sect.~5.1]{Sutti_V:2021}},
\begin{equation}\label{eq:var_pb_1}
   \begin{cases}
      \displaystyle\min_{w} \cF(w(x,y)) = \int_{\Omega} \tfrac{1}{2} \| \nabla w(x,y) \|^{2} - \gamma(x,y)\,w(x,y) \dx\dy \\
      \quad \text{such that} \quad w=0 \ \text{on} \ \partial\Omega,
   \end{cases}
\end{equation}
where $ \nabla = \big( \frac{\partial}{\partial x}, \frac{\partial}{\partial y} \big) $, $ \Omega = [0,1]^{2} $ and $ \gamma $ is the source term. 
The variational derivative (Euclidean gradient) of $ \cF $ is
\begin{equation}\label{eq:var_pb_1_grad}
   \frac{\delta \cF}{\delta w} = -\Delta w - \gamma.
\end{equation}
A critical point of~\eqref{eq:var_pb_1} is thus also a solution of the elliptic PDE $ -\Delta w = \gamma$.
We refer the reader to~\protect{\cite[Sect.~5.1.1]{Sutti_V:2021}} or \protect{\cite[Sect.~7.4.1.1]{Sutti:2020b}} for the details about the discretization.

Finally, we consider the variational problem from \protect{\cite[Sect.~5.2]{Sutti_V:2021}}:
\begin{equation}\label{eq:modified_NPDE_variational}
   \begin{cases}
      \displaystyle\min_{w} \cF(w) = \displaystyle\int_{\Omega} \tfrac{1}{2} \| \nabla w \|^{2} + \lambda w^{2} \big( \tfrac{1}{3}w + \tfrac{1}{2} \big) - \gamma \, w \dx\dy \\
      \quad \text{such that} \quad w=0 \ \text{on} \ \partial\Omega.
   \end{cases}
\end{equation}
For $ \gamma $, we choose
\[
   \gamma(x,y) = e^{x-2y} \, \sum_{j=1}^{5} 2^{j-1} \sin( j\pi x) \sin(j \pi y).
\]
which is the same right-hand side adopted in~\protect{\cite{Sutti_V:2021}}. The variational derivative of $ \cF $ is
\[
   \frac{\delta \cF}{\delta w} =  -\Delta w + \lambda w (w+1) - \gamma  = 0.
\]
Regardless of the specific form of the functional $ \cF $, all the problems studied in this paper have the general formulation
\[
   \min_{W} F(W) \ \text{s.t.} \ W \in \lbrace X \in \Rnn \colon \rank(X) = r \rbrace,
\]
where $ F $ denotes the discretization of the functional $ \cF $.

More details about each problem are provided later in the dedicated sections.

\section{Riemannian optimization framework} \label{sec:riem_optim}
As anticipated above, in this paper, we use the \textit{Riemannian optimization framework} \protect{\cite{Edelman:1998,AMS:2008}}. 
This approach exploits the underlying geometric structure of the low-rank constrained problems, thereby allowing to take explicitly into account the constraints.
In practice, the optimization variables in our variational problems are constrained to a smooth manifold, and we perform the optimization on the manifold.

Specifically, in this paper, we use the Riemannian trust-region (RTR) method of~\protect{\cite{ABG:2007}}. A more recent presentation of the RTR method can be found in~\protect{\cite{boumal_2023}}. In the next section, we introduce some fundamental geometry concepts used in Riemannian optimization, which are needed to formulate the RTR method, whose pseudocode is provided in Sect.~\ref{sec:RTR}.

\section{Geometry of the manifold of fixed-rank matrices} \label{sec:geom_fixed_rank_manif}
The manifold of fixed-rank matrices is defined as
\[
    \cMr = \lbrace X \in \Rmn \colon \rank(X) = r \rbrace.
\]
Using the singular value decomposition (SVD), one has the equivalent characterization 
\begin{equation*}
    \begin{split}
       \cMr = \lbrace U\Sigma V\tr \colon & U \in \Stmr, \ V \in \Stnr, \\
              & \Sigma = \diag(\sigma_{1}, \sigma_{2}, \ldots, \sigma_{r} ) \in \R^{r \times r}, \ \sigma_{1} \geq \cdots \geq \sigma_{r} > 0 \rbrace,
    \end{split}
\end{equation*}
where $ \Stmr $ is the Stiefel manifold of $ m \times r $ real matrices with orthonormal columns, and  $ \diag(\sigma_{1}, \sigma_{2}, \ldots, \sigma_{r} ) $ is a square matrix with $ \sigma_{1}, \sigma_{2}, \ldots, \sigma_{r} $ on its diagonal.

\subsection{Tangent space and metric}
The following proposition shows that $ \cMr $ is a smooth manifold with a compact representation for its tangent space.

\begin{proposition}[\protect{\cite[Prop.~2.1]{Vandereycken:2013}}]
   The set $ \cMr $ is a smooth submanifold of dimension $ (m+n-r)r $ embedded in $ \Rmn $. Its tangent space $ \TXMr $ at $ X = U\Sigma V\tr \in \cMr $ is given by
   \begin{equation}\label{eq:tan_vec_format}
      \TXMr =
      \begin{bmatrix}
         U  &  U_{\perp}
      \end{bmatrix}
      \begin{bmatrix}
         \R^{r \times r}      &  \R^{r \times (n-r)}  \\
         \R^{(m-r) \times r}  &  0_{(m-r)\times (n-r)}
      \end{bmatrix}
      \begin{bmatrix}
         V  &  V_{\perp}
      \end{bmatrix}\tr.  
   \end{equation}
   In addition, every tangent vector $\xi \in \TXMr$ can be written as
\begin{equation}\label{eq:tan_vec_format_small_param}
\xi = UMV\tr + \Up V\tr + U\Vp\tr,
\end{equation}
with $M \in \R^{r \times r}$, $\Up \in \R^{m\times r}$, $\Vp \in \R^{n \times r}$ such that $\Up\tr U = \Vp\tr V = 0$.
\end{proposition}
Since $ \cMr \subset \Rmn $, we represent tangent vectors in~\eqref{eq:tan_vec_format} and~\eqref{eq:tan_vec_format_small_param} as matrices of the same dimensions.
The Riemannian metric is the restriction of the Euclidean metric on $ \Rmn $ to the submanifold $ \cMr $,
\[
    g_{X}(\xi,\eta) = \langle \xi, \eta \rangle = \trace(\xi\tr \eta), \quad \text{with} \ X \in \cMr \ \text{and} \ \xi, \eta \in \TXMr.
\]

\subsection{Projectors}
Defining $ P_{U} = UU\tr $ and $ P_{U}^{\perp} = I - P_{U} $ for any $ U \in \Stmr $, where $ \Stmr $ is the Stiefel manifold of $m$-by-$ r $ orthonormal matrices, the orthogonal projection onto the tangent space at $ X $ is \protect{\cite[(2.5)]{Vandereycken:2013}}
\begin{equation}\label{eq:def_proj}
   \PTXM \colon \Rmn \to \TXMr, \quad Z \mapsto P_{U} Z P_{V} + P_{U}^{\perp} Z P_{V} + P_{U} Z P_{V}^{\perp}.
\end{equation}
Since this projector is a linear operator, we can represent it as a matrix.
The projection matrix $\PXmat\in \R^{n^{2} \times n^{2}}$ representing the operator $\PTXM$ can be written as
\[
    \PXmat \coloneqq P_{V} \otimes P_{U} + P_{V} \otimes P_{U}^{\perp} + P_{V}^{\perp} \otimes P_{U}.
\]

\subsection{Riemannian gradient}
The Riemannian gradient of a smooth function $ f \colon \cMr \to \R $ at $ X \in \cMr $ is defined as the unique tangent vector $ \grad f(X) $ in $ \TXMr $ such that
\[
   \forall \xi \in \TXMr, \quad \langle \, \grad f(X), \xi \, \rangle = \D f(X) [\xi],
\]
where $ \D f $ denotes the directional derivatives of $ f $.
More concretely, for embedded submanifolds, the Riemannian gradient is given by the orthogonal projection onto the tangent space of the Euclidean gradient of $ f $ seen as a function on the embedding space $ \Rmn $; see, e.g., \protect{\cite[(3.37)]{AMS:2008}}.
Then, denoting $ \nabla f(X) $ the Euclidean gradient of $ f $ at $ X $, the Riemannian gradient is given by
\begin{equation}\label{eq:Riemannian_grad_as_projection}
    \grad f(X) = \PTXM \big( \nabla f(X) \big).
\end{equation}
In other terms, the Riemannian gradient is given by the orthogonal projection of the Euclidean gradient onto the tangent space.

\subsection{Riemannian Hessian}
The Riemannian Hessian is defined by (see, e.g., \protect{\cite[def.~5.5.1]{AMS:2008}}, \protect{\cite[def.~5.14]{boumal_2023}})
\[
   \Hessian f(x) [\xi_{x}] = \nabla_{\xi_{x}} \, \grad f(x),
\]
where $ \nabla_{\xi_{x}} $ is the Levi-Civita connection.
If $\cM$ is a Riemannian submanifold of the Euclidean space $ \R^{n} $, as it is the case for the manifold of fixed-rank matrices, it follows that \protect{\cite[cor.~5.16]{boumal_2023}}
\[
   \forall \xi \in \mathrm{T}_{x}\cM, \qquad \Hessian f(x) [\xi] = \PTXM\big( \D \grad f(x)[\xi]\big).
\]
In practice, this is what we use in the calculations.

\subsection{Retraction} \label{sec:retraction}
We need a mapping to map the updates in the tangent space onto the manifold. This mapping is provided by a retraction $ \Retraction_{X} $. A retraction is a smooth map from the tangent space to the manifold, $ \Retraction_{X}\colon \TXMr \to \cMr $, used to map tangent vectors to points on the manifold. It is, essentially, any smooth first-order approximation of the exponential map of the manifold; see, e.g., \protect{\cite{Absil:2012}}. To establish convergence of the Riemannian algorithms, it is sufficient for the retraction to be defined only locally.

An excellent survey on low-rank retractions is given in~\protect{\cite{Absil:2015}}. In our setting, we have chosen the metric projection, which is provided by a truncated SVD.

\section{The RTR method} \label{sec:RTR}
As we anticipated above, to solve the implicit equation resulting from the time-integration scheme, we employ the Riemannian trust-region method of~\protect{\cite{ABG:2007}}. For reference, we provide the pseudocode for RTR in Algorithm~\ref{algo:RTR}. Step 4 in Algorithm~\ref{algo:RTR} uses the truncated conjugate gradient (tCG) of~\protect{\cite{Toint:1981,Steihaug:1983}}. This method lends itself very well to being preconditioned.

\begin{algorithm}
    \caption{Riemannian trust-region}\label{algo:RTR}
    Given $ \bar{\Delta} > 0, \, \Delta_{1} \in (0,\bar{\Delta}) $
    
        \For{$ i = 1,2, \dots $}{
           \textbf{Define} the second-order model
           \[
              m_{i} \colon \mathrm{T}_{x_{i}} \cM \to \R, \ \xi \mapsto f(x_{i}) + \left\langle \grad f(x_{i}), \xi \right\rangle + \frac{1}{2} \left\langle \Hessian f(x_{i})[\xi], \xi \right\rangle.
           \]
           
        \textbf{TR subproblem}: compute $ \eta_{i} $ by solving
        \[
           \eta_{i} = \argmin m_{i}(\xi) \ \text{s.t.} \ \|\xi \| \leq \Delta_{i}.
        \]
        
        Compute $ \rho_{i} = (\widehat{f}(0) - \widehat{f}_{i}(\eta_{i}))/(m_{i}(0) - m_{i}(\eta_{i})) $.
        
  \uIf{$ \rho_{i} \geq 0.05 $}{
    Accept step and set $ x_{i+1} = \Retraction_{x_{i}}(\eta_{i}) $.
  }
  \Else{
    Reject step and set $ x_{i+1} = x_{i} $.
  }
  Radius update: set
  \[
     \Delta_{i+1} =
     \begin{cases}
        \min(2\Delta_{i}, \bar{\Delta}) & \text{if} \ \rho_{i} \geq 0.75 \ \text{and} \ \| \eta_{i} \| = \Delta_{i}, \\
        0.25 \, \| \eta_{i} \| & \text{if} \ \rho_{i} \leq 0.25, \\
        \Delta_{i} & \text{otherwise}.
     \end{cases}     
  \]
  }
\end{algorithm}

\subsection{Riemannian gradient and Riemannian Hessian}
In general, in the case of Riemannian submanifolds, the full Riemannian Hessian of an objective function $f$ at $ x \in \cM $ is given by the projected Euclidean Hessian plus the curvature part
\begin{equation}\label{eq:full_riem_hess}
   \Hessian f(x) [\xi] = P_{x} \, \nabla^{2} f(x) \, P_{x} + P_{x} \, (\text{``curvature terms''}) \, P_{x}.
\end{equation}
This suggests using $ P_{x} \, \nabla^{2} f(x) \, P_{x} $ as a preconditioner in the Riemannian trust-region scheme. 
For the LYAP problem, the Riemannian gradient is given by
\begin{equation*}
    \grad F(X)  = \PTXM \! \big( h_{x}^{2} \, (AX + XA - \Gamma ) \big).
\end{equation*}

The directional derivative of the gradient, i.e., the Euclidean Hessian applied to $ \xi \in \TXMr $, is
\[
   \Hessian F(X) [\xi] = \D \grad F(X)[\xi] = h_{x}^{2} \, (A \xi + \xi A).
\]
The orthogonal projection of the Euclidean Hessian followed by vectorization yields
\begin{align*}
   \vecop \big(\PTXM \! \big( \D \grad F(x)[\xi]\big)\big) &= h_{x}^{2} \, \PXmat \vecop (A \xi + \xi A) \\
                                    &= h_{x}^{2} \, \PXmat (A \otimes I + I \otimes A)  \vecop(\xi) \\
                                    &= h_{x}^{2} \, \PXmat (A \otimes I + I \otimes A) \, \PXmat \vecop(\xi),
\end{align*}
where the second $ \PXmat $ is inserted for symmetrization.
From here we can read the symmetric $n^{2}$-by-$n^{2}$ matrix
\begin{equation}\label{eq:prec_candidate}
   H_{X} = h_{x}^{2} \, \PXmat (A \otimes I + I \otimes A) \, \PXmat.
\end{equation}
The inverse of this matrix~\eqref{eq:prec_candidate} should be a good candidate for a preconditioner. In the next section, we present the derivation of the preconditioner on the manifold of fixed-rank matrices.

In general, the preconditioner from above cannot be efficiently inverted because of the coupling with the nonlinear terms. Nonetheless, numerical experiments in the following sections show that it remains an efficient preconditioner even for problems with a (mild) nonlinearity.

\section{The Allen--Cahn equation} \label{sec:ACE}

The Allen--Cahn equation is a reaction-diffusion equation originally studied for modeling the phase separation process in multi-component alloy systems~\protect{\cite{Allen:1972,Allen:1973}}. It later turned out that the Allen--Cahn equation has a much wider range of applications. Recently, \protect{\cite{Yoon:2020}} provided a good review. Applications include mean curvature flows~\protect{\cite{Laux:2018}}, two-phase incompressible fluids~\protect{\cite{YangFeng:2006}}, complex dynamics of dendritic growth~\protect{\cite{LeeKim:2012}}, image inpainting~\protect{\cite{Dobrosotskaya:2008,LiJeong:2015}}, and image segmentation~\protect{\cite{Benes:2004,Lee:2019}}.

The Allen--Cahn equation in its simplest form reads 
\begin{equation}\label{eq:allen-cahn}
   \frac{\partial w}{\partial t} = \varepsilon \Delta w + w  - w^{3},
\end{equation}
where $ w \equiv w(\bm{x}, t) $, $ \bm{x} \in \Omega = [-\pi, \pi]^{2} $, and $ t\geq 0 $. It is a stiff PDE with a low-order polynomial nonlinearity and a diffusion term $ \varepsilon \Delta w $. As in~\protect{\cite{RodgersVenturi:2022}}, we set $ \varepsilon = 0.1 $, and we solve \eqref{eq:allen-cahn} on a two-dimensional flat torus, and we also use the same initial condition as in~\protect{\cite[(77)--(78)]{RodgersVenturi:2022}}, namely,
\begin{equation}\label{eq:ACE_initial_condition}
   w_{0}(x,y) = u(x, y) - u(x, 2y) + u(3x + \pi, 3y + \pi) - 2 u(4x, 4y) + 2 u(5x, 5y).
\end{equation}
where
\[
   u(x,y) = \frac{\left[ e^{-\tan^{2}(x)} + e^{-\tan^{2}(y)}\right] \sin(x) \sin(y)}{1 + e^{\lvert \csc (-x/2)\rvert}  + e^{\lvert \csc (-y/2)\rvert}}.
\]
We emphasize that with this choice, the matrix $ W_{0} $ has no low-rank structure and will be treated as a dense matrix. Similarly, during the first 0.5 seconds of the simulation, the numerical solution $ W $ will be treated as a full-rank matrix. Nonetheless, thanks to the Laplacian's smoothing effect as time evolution progresses, the solution $ W $ can be well approximated by low-rank matrices. In particular, for large simulation times, the solution converges to either 0 or 1 in most of the domain (cf. panel (f) of Figure~\ref{fig:123456}).

\subsection{Spatial discretization}
We discretize~\eqref{eq:allen-cahn} in space on a uniform grid with $ 256 \times 256 $ points. 
In particular, we use the central finite differences to discretize the Laplacian with periodic boundary conditions.
This results in the matrix ODE
\begin{equation}\label{eq:allen-cahn_discretized}
   \dt{W} = \varepsilon \left( A W + W  A \right) + W - W^{\circ 3},
\end{equation}
where $W \colon [0, T] \to \R^{256 \times 256} $ is a matrix that depends on $ t $, $^{\circ 3}$ denotes the elementwise power of a matrix (so-called Hadamard power, defined by $ W^{\circ \alpha} = [w_{ij}^{\alpha}] $), and $ A $ is the second-order periodic finite difference differentiation matrix
\begin{equation}\label{eq:discretized_laplacian}
   A = \frac{1}{h_{x}^{2}}
   \begin{bmatrix}
      -2  &       1  &          &         &    1 \\
       1  &      -2  &       1  &         &      \\
          &  \ddots  &  \ddots  & \ddots  &      \\
          &          &       1  &     -2  &    1 \\
       1  &          &          &      1  &   -2
   \end{bmatrix}.
\end{equation}

This matrix ODE is an initial value problem (IVP) in the form of \protect{\cite[(48)]{UschmajewV:2020}}
\begin{equation}\label{eq:IVP}
   \begin{cases}
      \dt{W}(t) = G(W(t)), \\
      W(t_{0}) = W_{0},
   \end{cases}
\end{equation}
where $ G \coloneqq \varepsilon \left( A W + W  A \right) + W - W^{\circ 3} $ is the right-hand side of~\eqref{eq:allen-cahn_discretized}, and $ W_{0} $ the discretization of the initial condition \eqref{eq:ACE_initial_condition}.

\subsection{Reference solution}

\begin{figure}[htbp]
    \centering
    \begin{minipage}{0.33\textwidth}
        \centering
        \includegraphics[width=\textwidth]{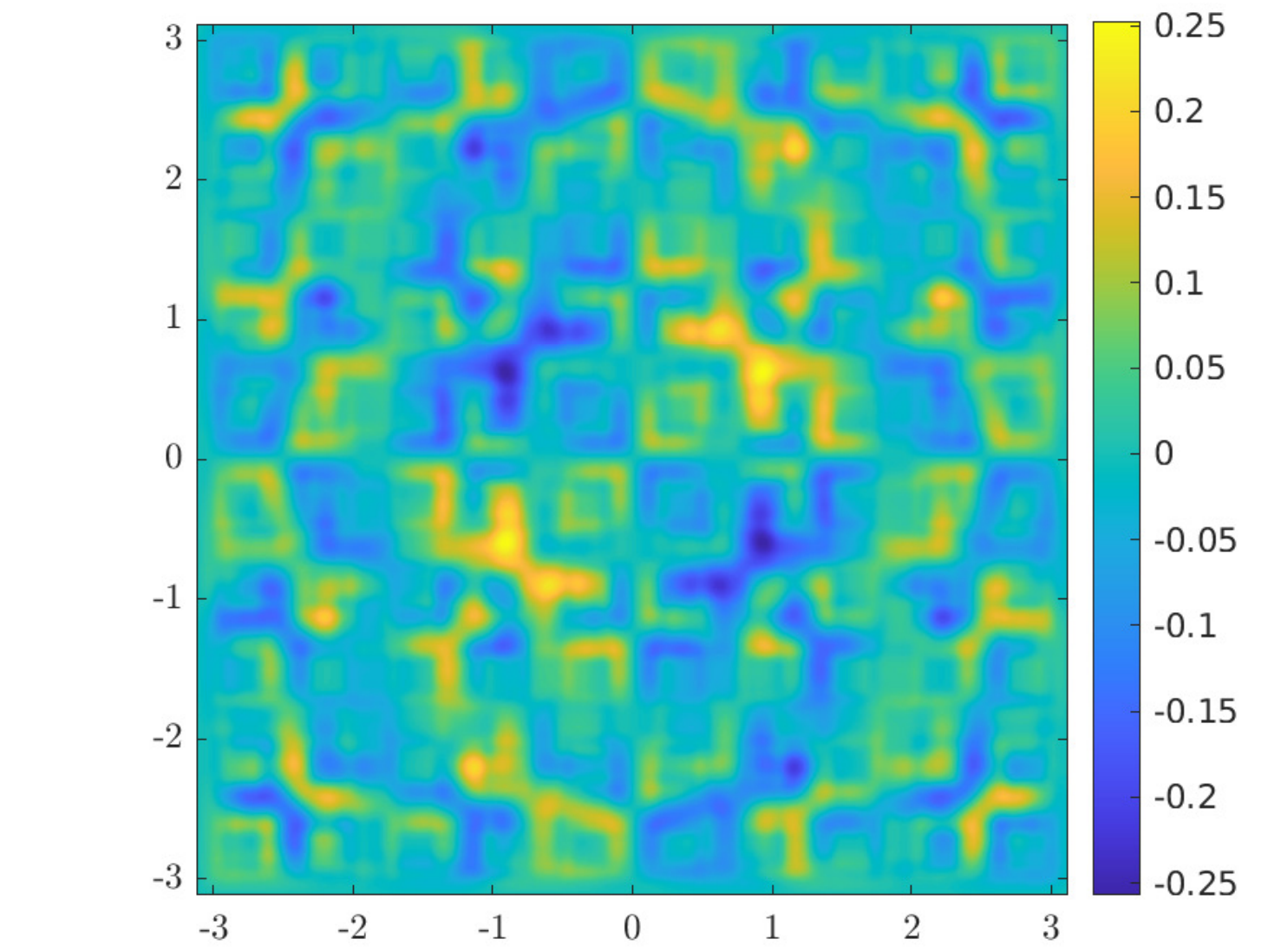}
        {\scriptsize (a) $ t = 0 $}
    \end{minipage}\hfill
    \begin{minipage}{0.33\textwidth}
        \centering
        \includegraphics[width=\textwidth]{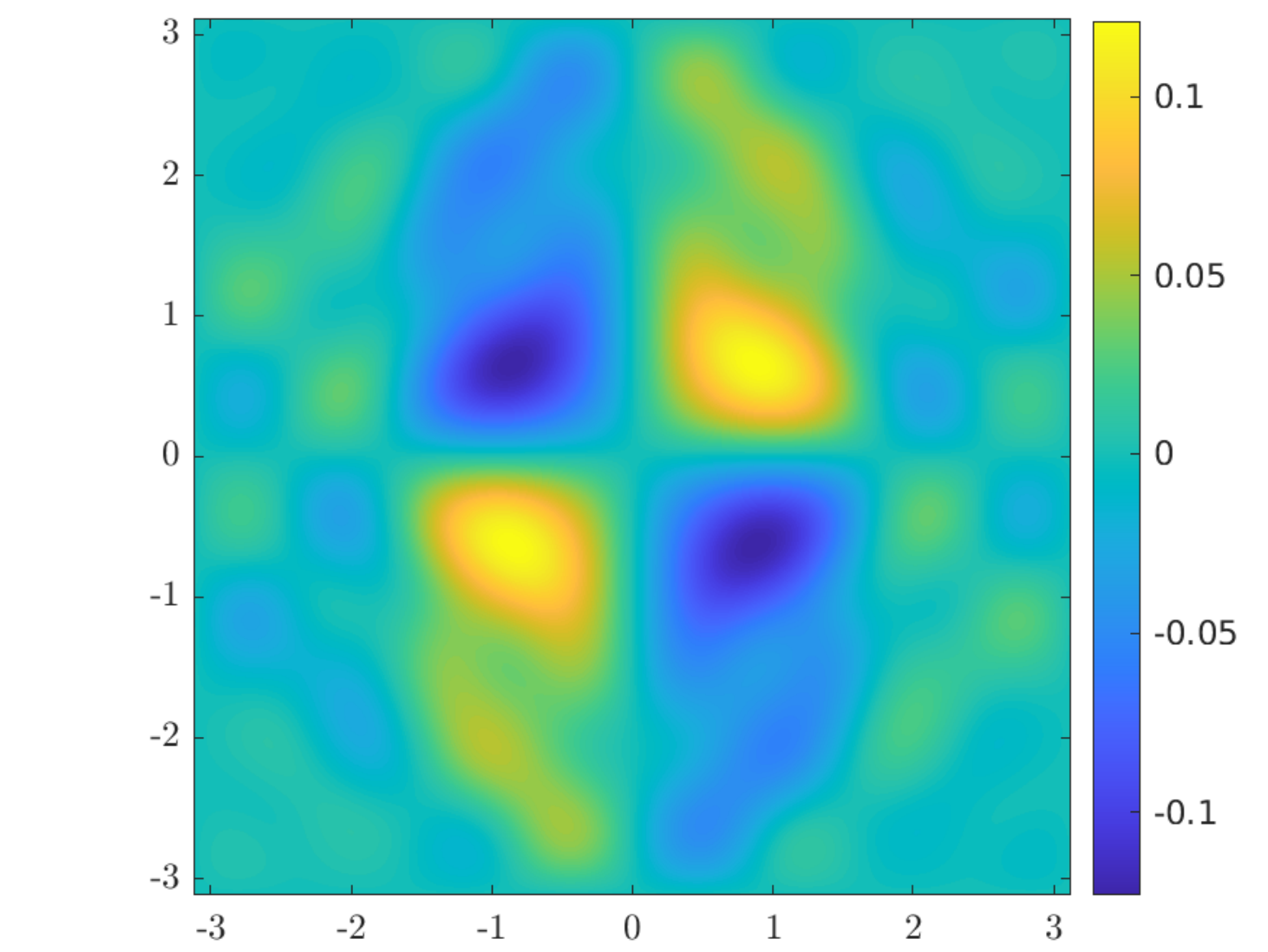}
        {\scriptsize (b) $ t = 0.5 $}        
    \end{minipage}\hfill
    \begin{minipage}{0.33\textwidth}
        \centering
        \includegraphics[width=\textwidth]{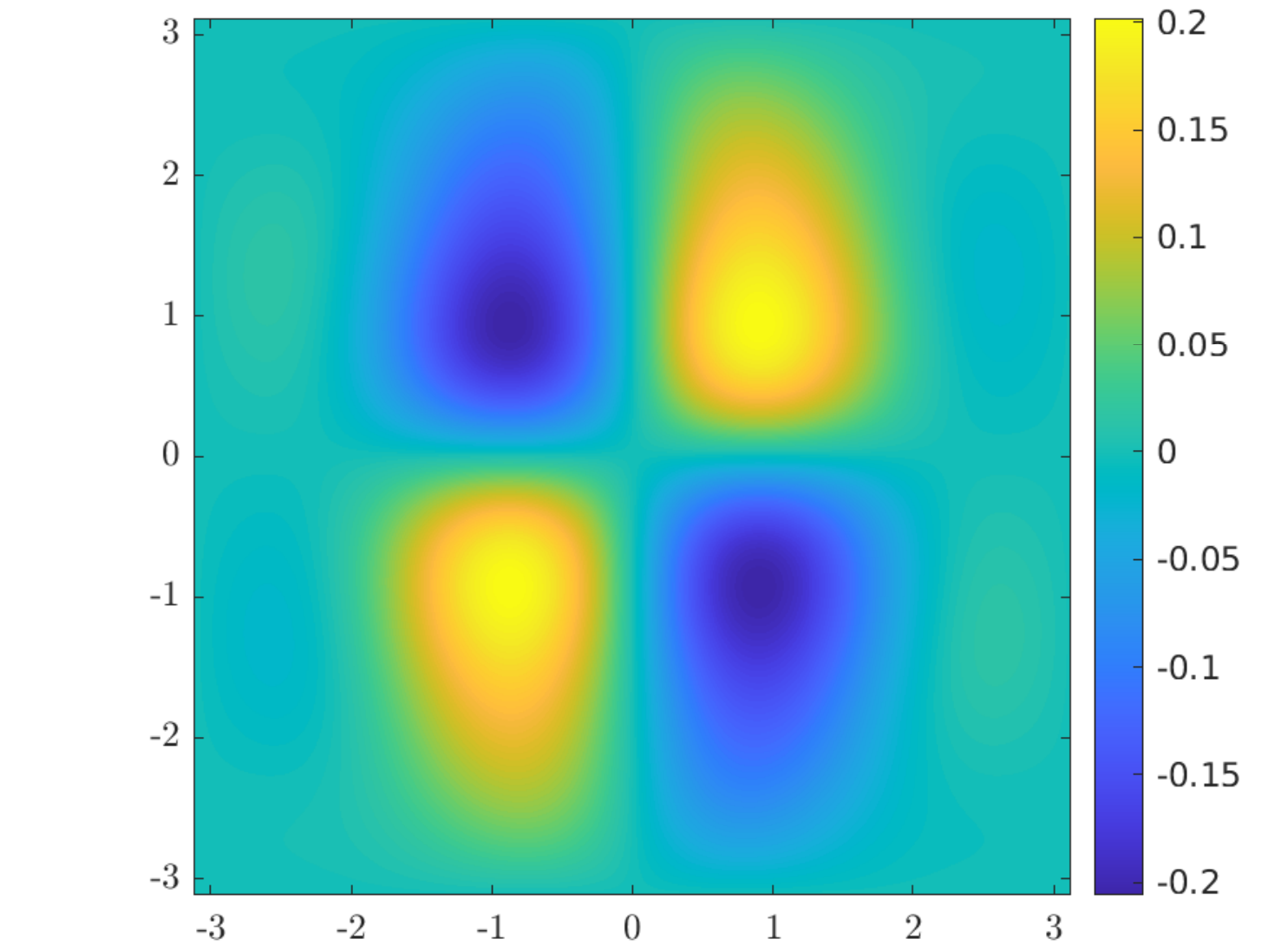}
        {\scriptsize (c) $ t = 2 $}        
    \end{minipage}\\
    \bigskip
    \centering
    \begin{minipage}{0.33\textwidth}
        \centering
        \includegraphics[width=\textwidth]{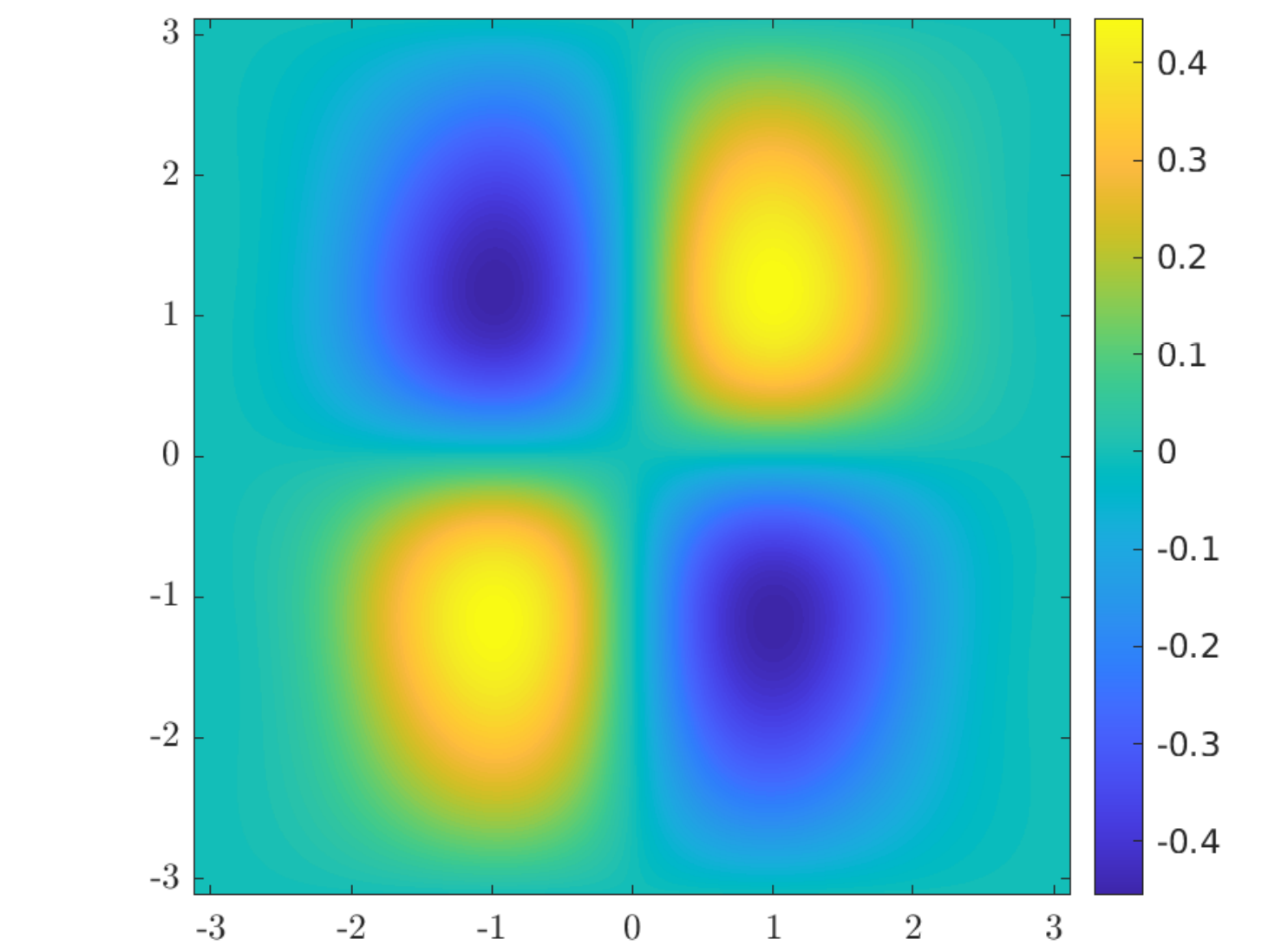}
        {\scriptsize (d) $ t = 3.5 $}        
    \end{minipage}\hfill
    \begin{minipage}{0.33\textwidth}
        \centering
        \includegraphics[width=\textwidth]{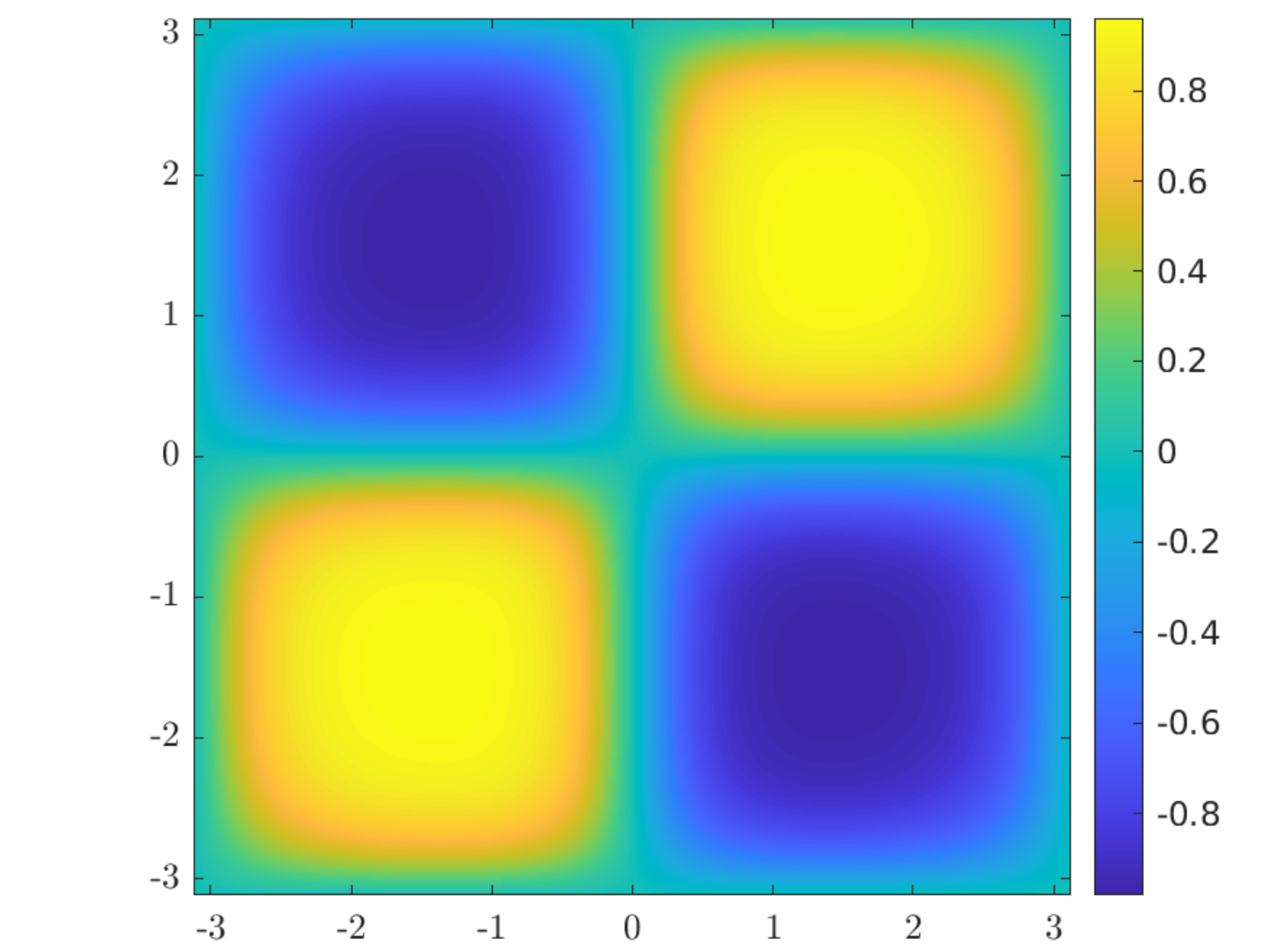}
        {\scriptsize (e) $ t = 7.5 $}        
    \end{minipage}\hfill
    \begin{minipage}{0.33\textwidth}
        \centering
        \includegraphics[width=\textwidth]{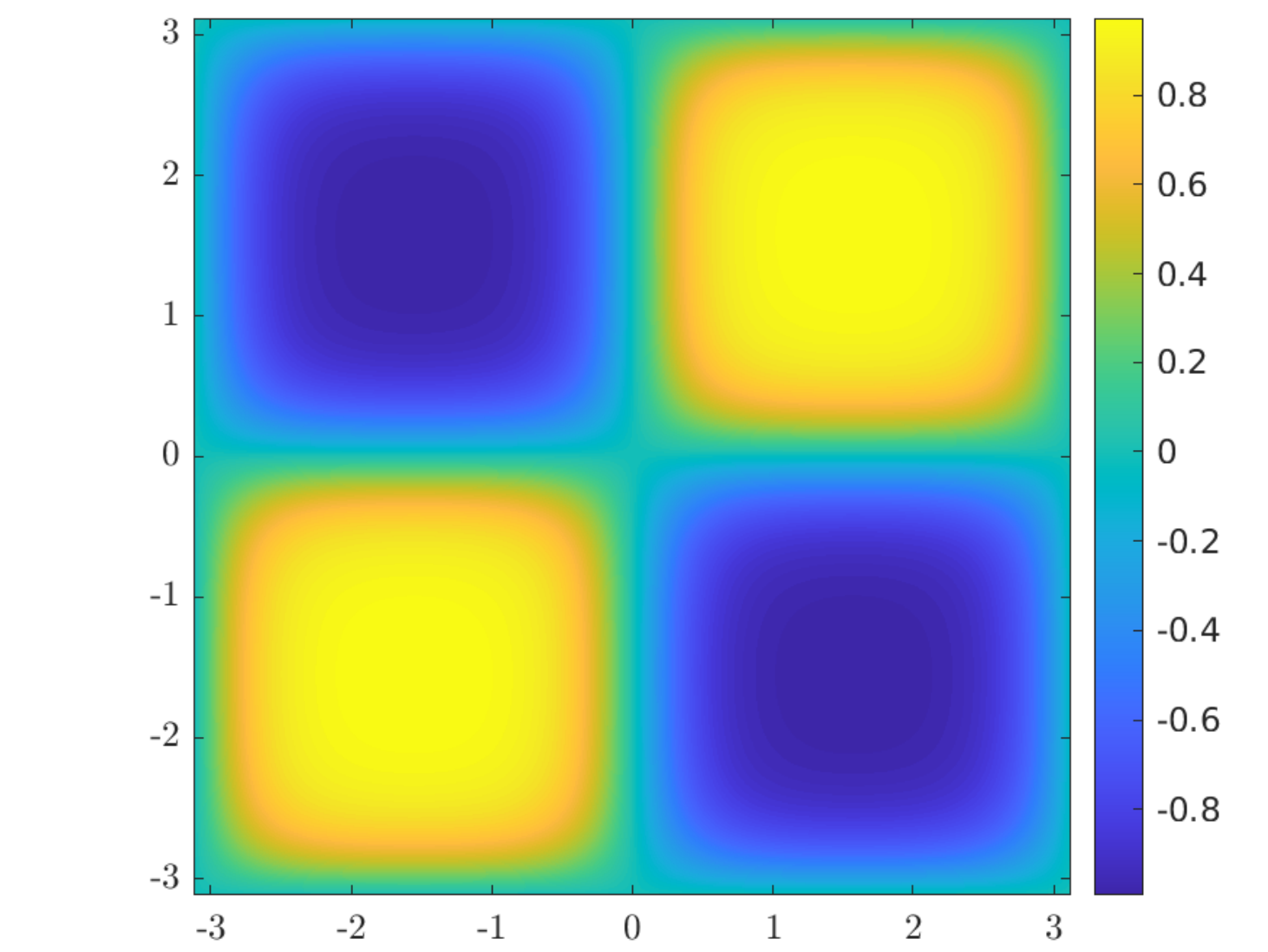}
        {\scriptsize (f) $ t = 15 $}
    \end{minipage}
    \caption{Time evolution of the solution $w$ to the Allen--Cahn equation, computed with ERK4, $ \deltatime = 10^{-4} $.}\label{fig:123456}
\end{figure}

To get a reference solution $ W_{\mathrm{ref}} $, we solve the (full-rank) IVP problem~\eqref{eq:IVP} with the classical explicit fourth-order Runge--Kutta method (ERK4), with a time step $ \deltatime = 10^{-4} $.
Figure~\ref{fig:123456} illustrates the time evolution of the solution to the Allen--Cahn equation.

\begin{figure}[htbp]
    \centering
    \begin{minipage}[c]{0.485\textwidth}
        \centering
        \includegraphics[width=0.97\textwidth]{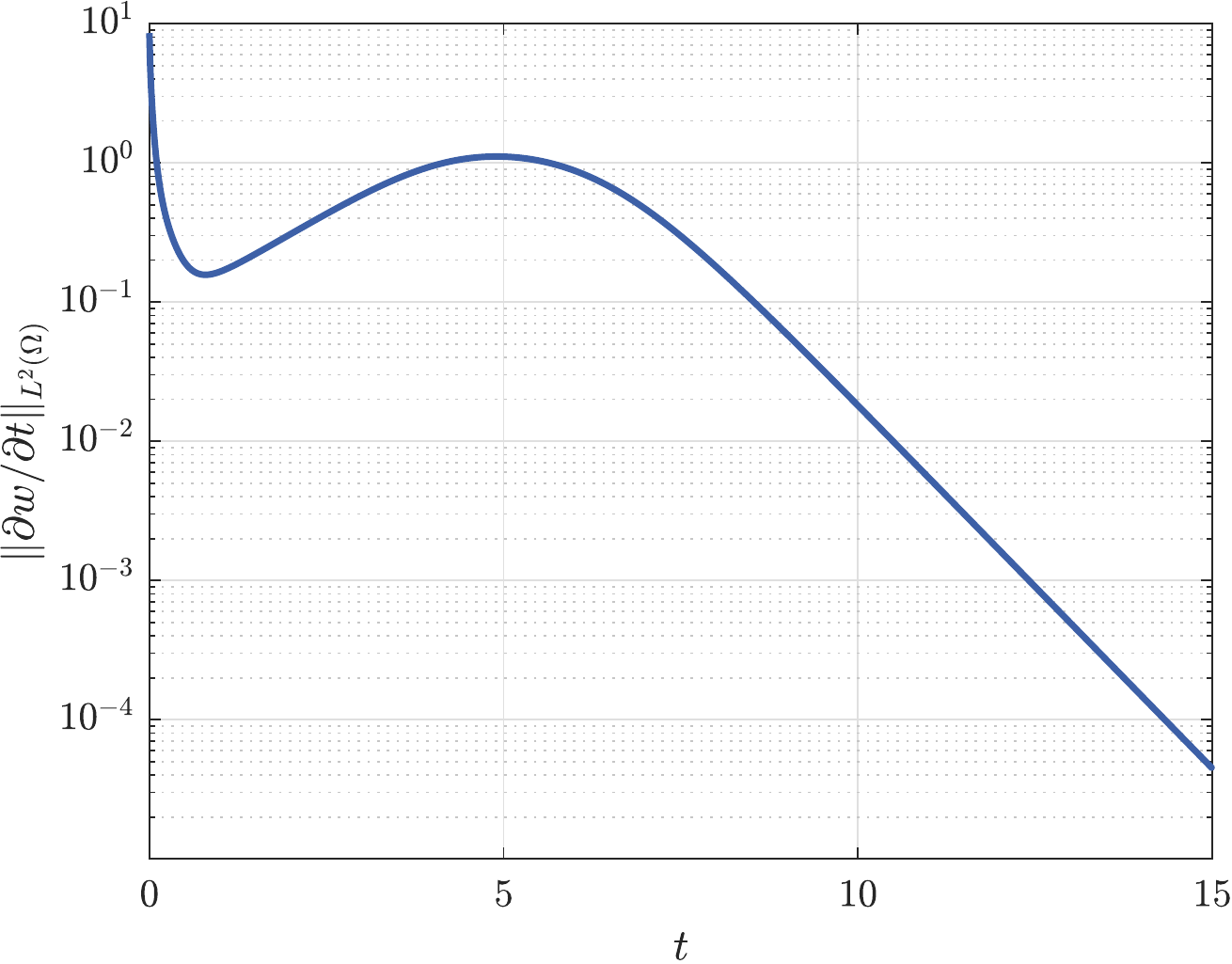}
        {\scriptsize (a) \emph{Discrete $ L^{2} $-norm of the RHS of~\eqref{eq:allen-cahn_discretized}}.}
    \end{minipage}\hfill
    \begin{minipage}[c]{0.485\textwidth}
        \centering
        \includegraphics[width=\textwidth]{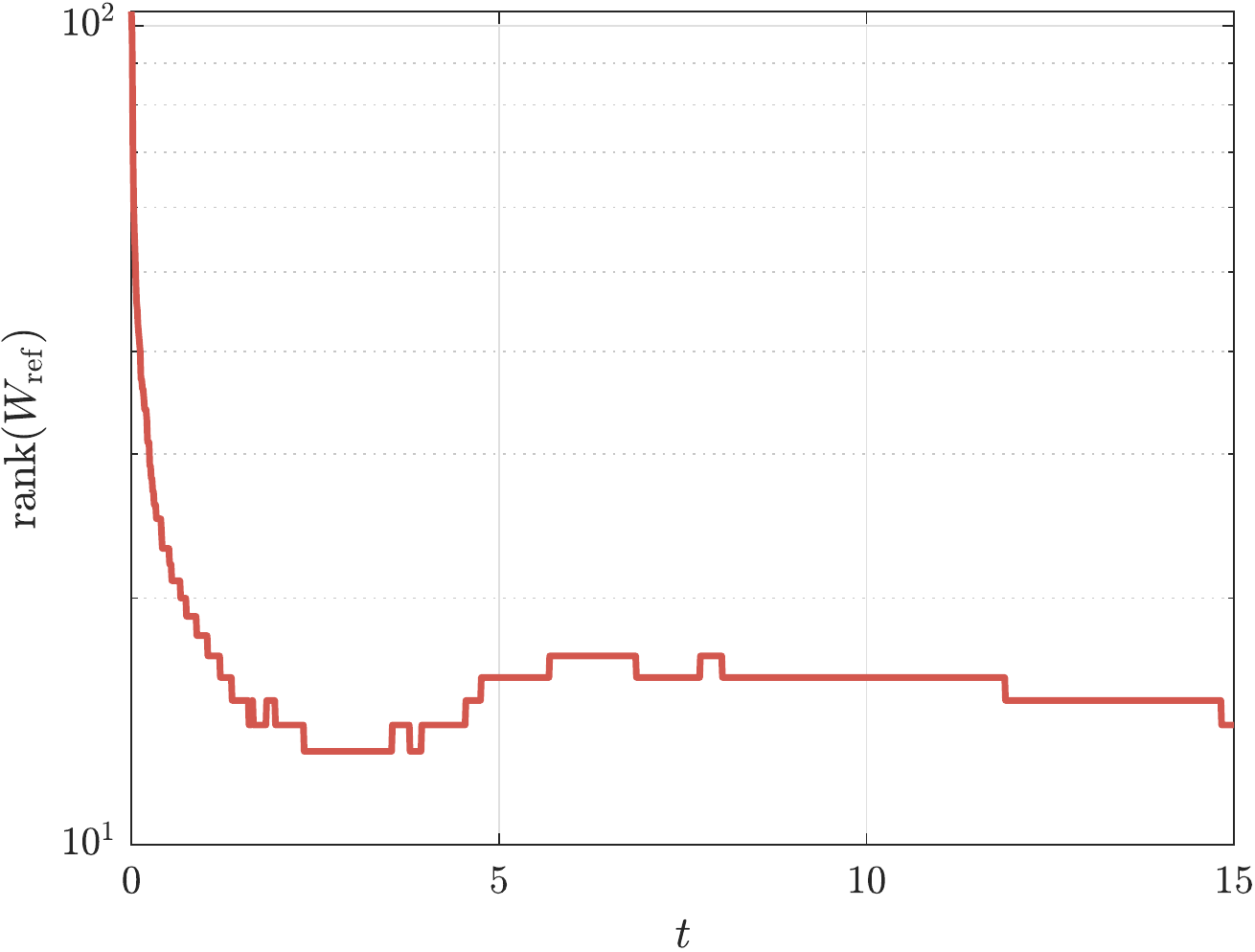}
        {\scriptsize (b) \emph{Numerical rank of $ W_{\mathrm{ref}} $}.}
    \end{minipage}
   \caption{Preliminary numerical study of $ W_{\mathrm{ref}} $.}\label{fig:78}
\end{figure}

As a preliminary study, we monitor the discrete $ L^{2} $-norm of the right-hand side of \eqref{eq:allen-cahn} for this reference solution. From panel (a) of Figure~\ref{fig:78}, it appears that after $ t \approx 13 $, $ \| \partial w / \partial t \|_{L^{2}(\Omega)} \approx 10^{-3} $, which means that the solution $ w $ enters a stationary phase (see also last two panels of Figure~\ref{fig:123456}).
Panel (b) of Figure~\ref{fig:78} plots the numerical rank of $ W_{\mathrm{ref}} $ versus time, with relative singular value tolerance of $ 10^{-10} $. The numerical rank exhibits a rapid decay during the first $ \approx 2 $ seconds, then varies between 13 and 17 during the rest of the simulation. The rank decreases as the diffusion term comes to dominate the system.

\subsection{Implicit time-stepping scheme}

As mentioned above, we employ the implicit Euler method for the time integration of~\eqref{eq:allen-cahn_discretized}, which gives
\begin{equation}\label{eq:implicit_euler}
   W_{k+1} = W_{k} + \deltatime \cdot G(W_{k+1}),
\end{equation}
and, additionally, we want $ W_{k} $ to be of low rank. This is done by using our PrecRTR on the manifold of fixed-rank matrices to solve for $ W_{k+1} $ the nonlinear equation~\eqref{eq:implicit_euler}.

Since our strategy is optimization, and we wish to maintain a kind of coherence with the problems LYAP and NPDE presented in the previous sections, we build a variational problem whose first-order optimality condition will be exactly~\eqref{eq:implicit_euler}. This leads us to consider the problem
\begin{equation}\label{eq:ACE_var_pb}
   \min_{w} \cF(w) \coloneqq \int_{\Omega} \frac{\varepsilon \deltatime}{2} \| \nabla w \|^{2} + \frac{(1-\deltatime)}{2} \, w^{2} + \frac{\deltatime}{4} \, w^{4} - \widetilde{w} \cdot w \dx \dy.
\end{equation}
It is interesting to note that this cost function is very similar to the NPDE functional \protect{\cite[(5.11)]{Sutti_V:2021}}. Here, $ \widetilde{w} $ is the solution at the previous time step and plays a similar role as $ \gamma $ in the NPDE functional (it is constant w.r.t. $ w $). Only the term $ \deltatime / 4 \, w^{4} $ is kind of novel w.r.t. NPDE. Moreover, in contrast to LYAP and NPDE, we need to solve this optimization problem many times, i.e., at every time step, to describe the time evolution of $ w $.
Our algorithm is summarized in Algorithm~\ref{algo:LRIE_ACE}.

\begin{algorithm}
  \caption{Low-rank Riemannian implicit Euler for the Allen--Cahn equation}\label{algo:LRIE_ACE}
   Given the initial condition~\eqref{eq:ACE_initial_condition}, $ \varepsilon $, and $ \deltatime $ and the RTR parameters\;
   $ k \leftarrow 0 $\;
   \While{until $ T $ is reached}{
      Solve (the discretized form of) \eqref{eq:ACE_var_pb} for $ W_{k} $ with PrecRTR on the manifold of fixed-rank matrices\;
      Set $ \widetilde{W} \leftarrow W_{k} $\;
   }
\end{algorithm}

We aim to obtain good low-rank approximations on the whole interval $ [t_{0}, T] $. However, it is clear from our preliminary study (panel (b) of Figure~\ref{fig:78}) that at the beginning of the time evolution, the numerical solution is not really low rank due to the initial condition of \protect{\cite{RodgersVenturi:2022}}. For this reason, in our numerical experiments, we consider the dense matrix until $ t = 0.5 $, and only then do we start our rank-adaptive method.

The discretizations of the objective function $ \cF(w) $ and its gradient are detailed in Appendix~\ref{app:discretization_ACE}.

\subsection{Numerical experiments}
The algorithm was implemented in MATLAB and will be publicly available at \url{https://github.com/MarcoSutti/PrecRTR}.
The Riemannian trust-region method of~\protect{\cite{ABG:2007}} was executed using solvers from the Manopt package~\cite{Boumal:2014} with the Riemannian embedded submanifold geometry from~\protect{\cite{Vandereycken:2013}}.
We conducted our experiments on a desktop machine with Ubuntu 22.04.1 LTS and MATLAB R2022a installed, with Intel Core i7-8700 CPU, 16GB RAM, and Mesa Intel UHD Graphics 630.

\begin{figure}[htbp]
    \centering
    \begin{minipage}[b]{0.485\textwidth}
        \centering
        \includegraphics[width=\textwidth]{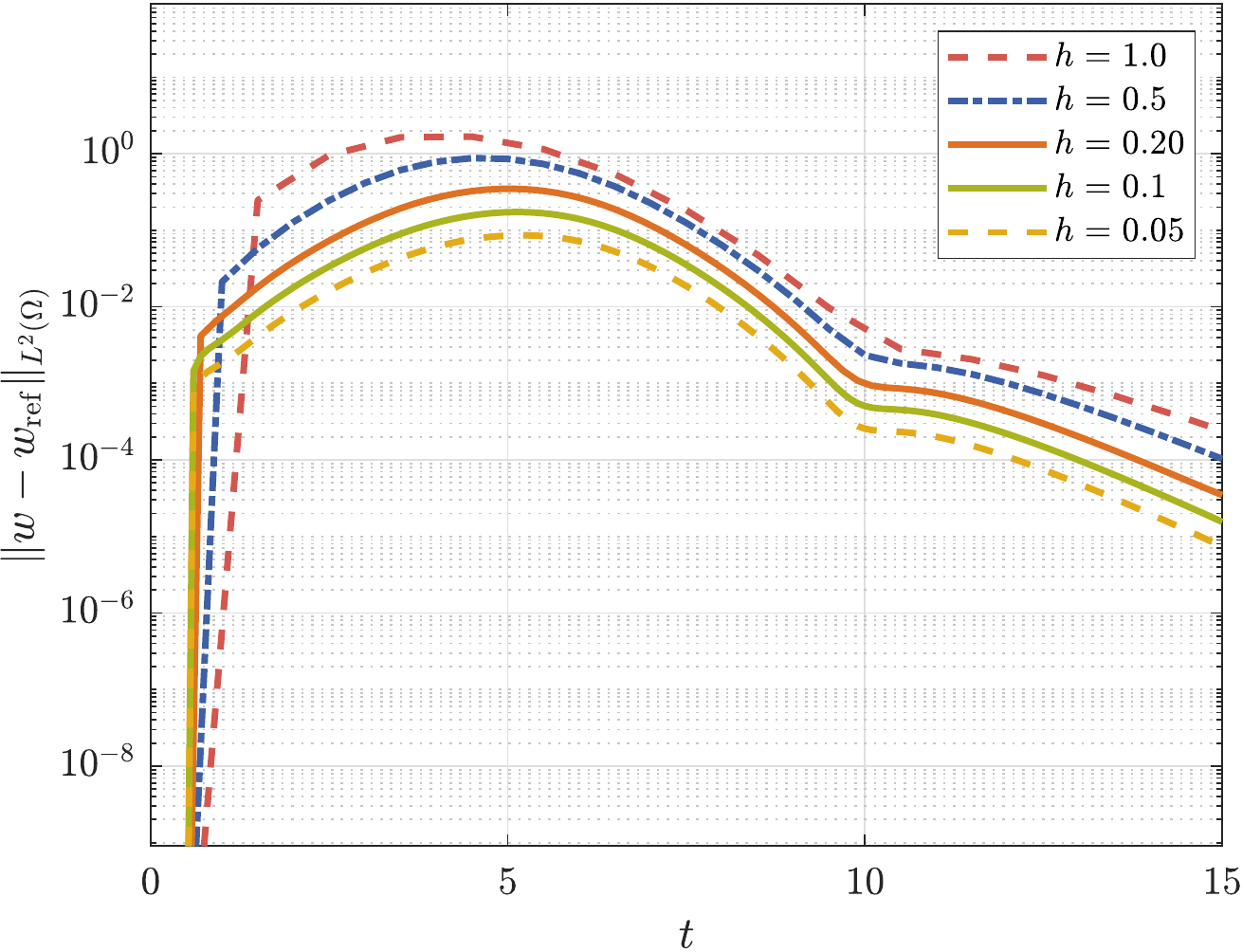}
        {\scriptsize (a)}
    \end{minipage}\hfill
    \begin{minipage}[b]{0.490\textwidth}
        \centering
        \includegraphics[width=\textwidth]{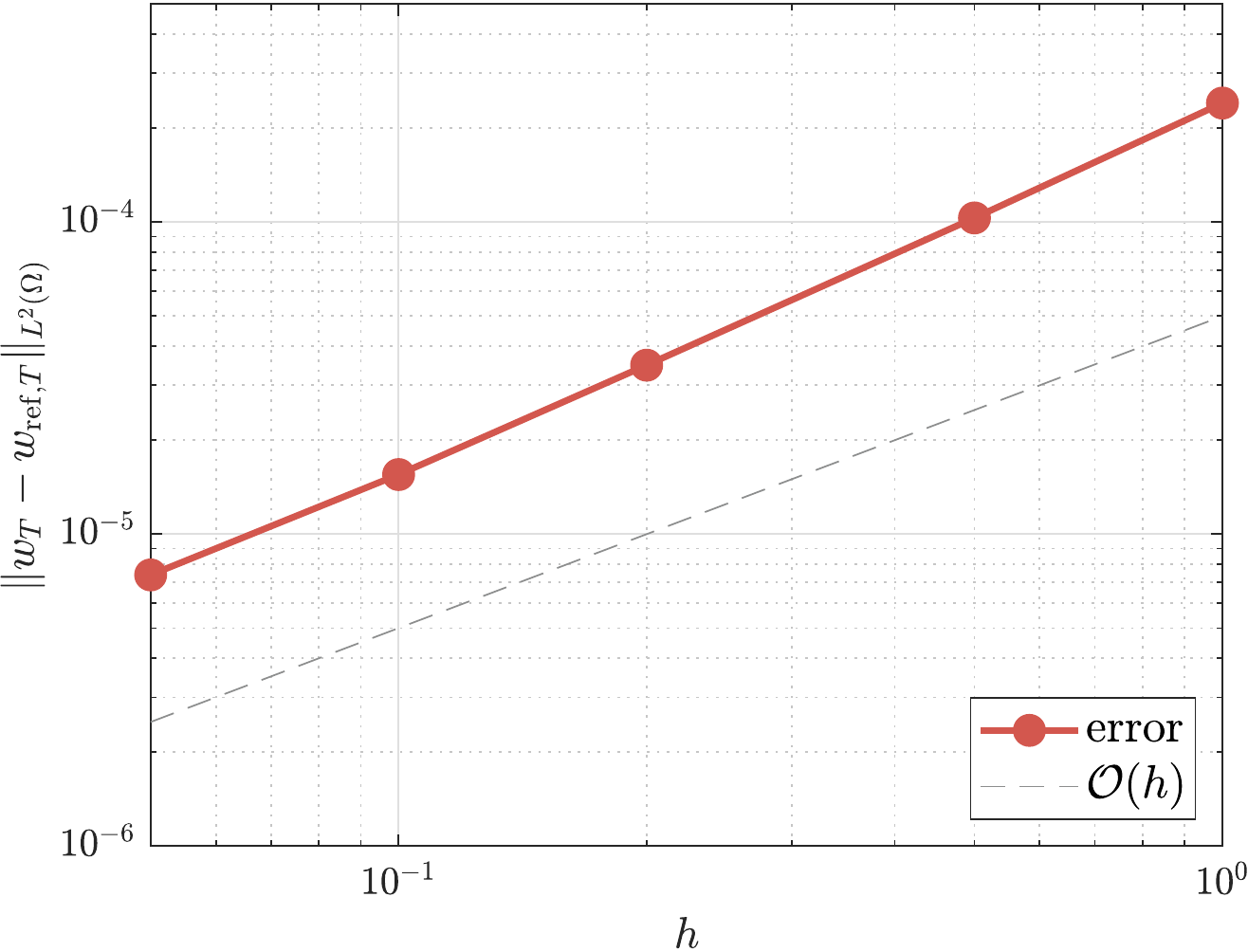}
        {\scriptsize (b)}
    \end{minipage}
    \caption{Panel (a): error versus time for the preconditioned low-rank evolution of the Allen-Cahn equation~\eqref{eq:allen-cahn} with initial condition~\eqref{eq:ACE_initial_condition}. Panel (b): error at $ T = 15 $ versus time step $ \deltatime $.}\label{fig:910}
\end{figure}

For the time integration, we use the time steps $ \deltatime = \lbrace 0.05, 0.1, 0.2, 0.5, 1 \rbrace $, and we monitor the error $ \| w - w_{\mathrm{ref}} \|_{L^{2}(\Omega)} $. 
In all expressions below, ${\cdot}^{(i)}$ indicates that a quantity was evaluated at the $i$th outer iteration of the Riemannian trust-region method.
Figure~\ref{fig:910} reports on the results.
Panel (a) shows the time evolution of the error $ \| w - w_{\mathrm{ref}} \|_{L^{2}(\Omega)} $, while panel (b) shows that the error decays linearly in $ \deltatime $.

\begin{figure}[htbp]
   \centering
   \includegraphics[width=0.5\columnwidth]{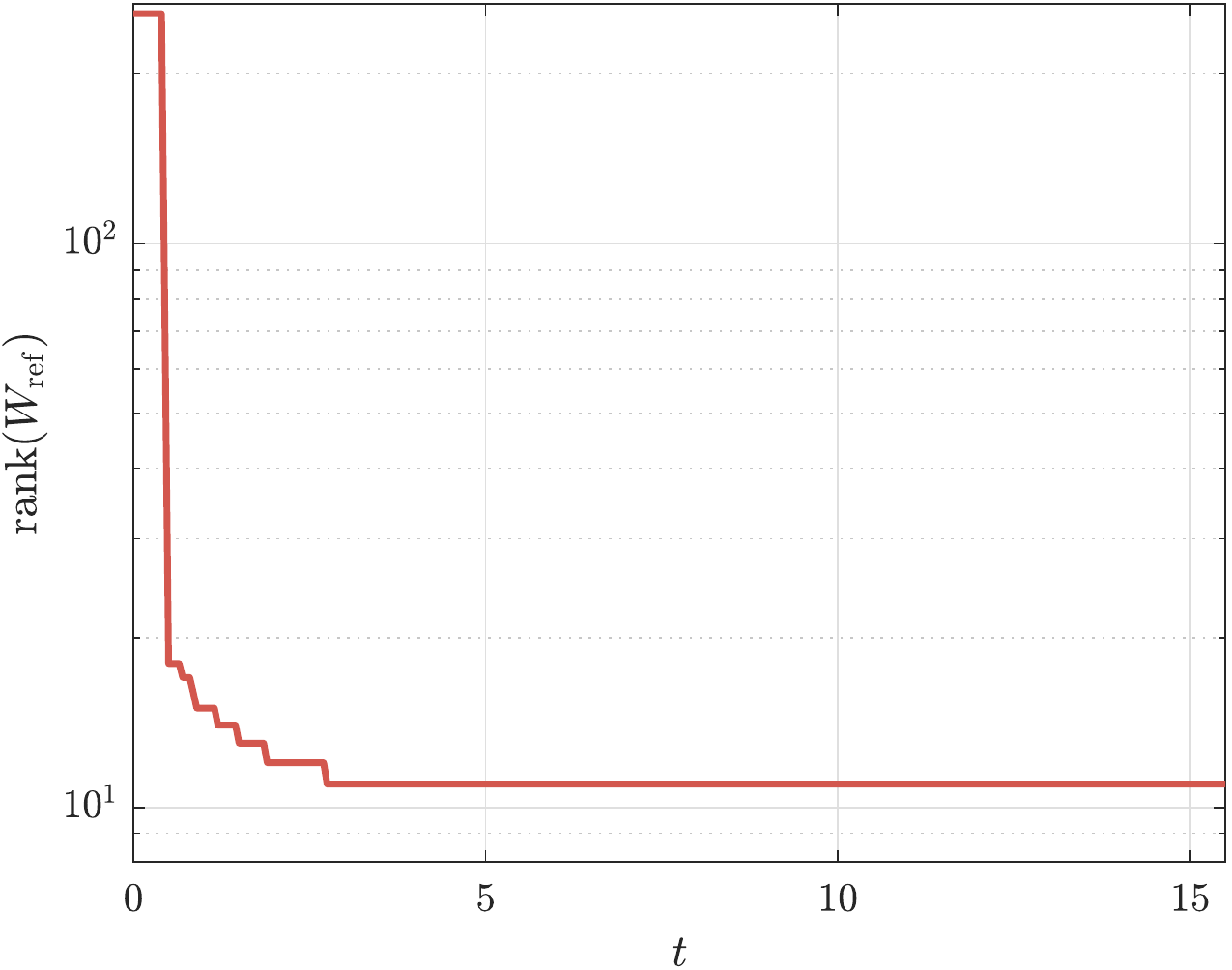}
   \caption{Rank versus time for the preconditioned low-rank evolution of Allen-Cahn equation~\eqref{eq:allen-cahn} with initial condition~\eqref{eq:ACE_initial_condition}, with $ \deltatime = 0.05 $.}\label{fig:11}
\end{figure}

Figure~\ref{fig:11} reports on the behavior of the rank with respect to time for the simulation with $ \deltatime = 0.05 $.

\subsection{Discussion/Comparison with other solvers}

It is evident that, even with very large time steps, we can still obtain relatively good low-rank approximations of the solution, especially at the final time $ T = 15 $. For example, compare with Figure 4 in~\protect{\cite{RodgersVenturi:2022}}, where the biggest time step is $ \deltatime = 0.01 $ --- i.e., one hundred times smaller than our largest time step. Moreover, in \protect{\cite{RodgersVenturi:2022}}, factorized formats are not mentioned. In contrast, we always work with the factors to reduce computational costs.

In~\protect{\cite{RodgersVenturi:2022}}, the authors study implicit rank-adaptive algorithms based on performing one time step with a conventional time-stepping scheme, followed by an implicit fixed-point iteration step involving a rank truncation operation onto a tensor or matrix manifold. Here, we also employ an implicit time-stepping scheme for the time evolution. Still, instead of using a fixed-point iteration method for solving the nonlinear equations, we use our preconditioned Riemannian trust regions (PrecRTR) on the manifold of fixed-rank matrices. 
This way, obtain a preconditioned dynamical low-rank approximation of the Allen--Cahn equation.

Implicit methods are much more effective for stiff problems.

Although implicit time integration methods are more expensive, since solutions of nonlinear systems replace function evaluations, they allow for a larger time step than their explicit counterparts. The additional computational overhead of the implicit method is compensated by the fact that we can afford a larger time step, as demonstrated by Figure~\ref{fig:910}.

The problem with explicit methods is that the time step needs to be in proportion to $ \sigma_{\mathrm{min}}(W(t)) $~\protect{\cite{UschmajewV:2020}}.

Moreover, the cost of solving the inner nonlinear equations remains moderately low thanks to our preconditioner.

When using fixed point iterations, one still obtains a condition on the time step size, which depends on the Lipschitz constant of the right-hand side term (see, e.g.,~\protect{\cite[(3.13)]{Kressner:2015}} and the reference cited therein~ \protect{\cite{Dieudonne:1960}})
\[
   h < \frac{1}{\| A \|_{\infty} \| G \|_{\infty}},
\]
where $ A $ is the matrix of the coefficients defining the stages (the Butcher tableau).
This condition appears to be a restriction on the time step, not better than the restrictions for explicit methods to be stable. This shows that (quote from~\protect{\cite{Kressner:2015}}) ``fixed point iterations are unsuitable for solving the nonlinear system defining the stages. For solving the nonlinear system, other methods like the Newton method should be used''.
A similar condition also holds for the method of Rodgers and Venturi, see~\protect{\cite[(31)]{RodgersVenturi:2022}}: 
\[
   h < \frac{1}{L_{G}}.
\]
Their paper states: ``Equation (31) can be seen as a stability condition restricting the maximum allowable time step $ \deltatime $ for the implicit Euler method with fixed point iterations.'' This makes a case for using the Newton method instead of fixed point iteration to find a solution to the nonlinear equation.

As we observed from the MATLAB profiler, as $ \ell $ increases, the calculation of the preconditioner becomes dominant in the running time. We are in the best possible situation since the preconditioner dominates the cost.

Low-rank Lyapunov solvers (see~\cite{Simoncini:2016} for a review) cannot be used to solve this kind of problem due to a nonlinear term in the Hessian, and PrecRTR proves much more effective than the RMGLS method of~\protect{\cite{Sutti_V:2021}}. However, the latter may remain effective in all those problems for which an effective preconditioner is unavailable.

Our method proves efficient when the rank is low and the time step is not too small. Otherwise, if these conditions are not met, there is no advantage over using full-rank matrices.

\section{The Fisher--KPP equation} \label{sec:FKPP}

The Fisher--KPP equation is a nonlinear reaction-diffusion PDE, which in its simplest form reads \protect{\cite[(13.4)]{Murray:2002}}
\begin{equation}\label{eq:FKPP_equation}
   \frac{\partial w}{\partial t} = \frac{\partial^{2} w}{\partial x^{2}} + r(\omega) \, w (1-w),
\end{equation}
where $ w \equiv w(x, t; \omega ) $, $ r( \omega ) $ is a species's reaction rate or growth rate. It is called ``stochastic''\footnote{``Stochastic'' might be too big of a term since no Brownian motion is involved. It is just a PDE with random coefficients for the initial condition and the reaction rate. } Fisher--KPP equation in the recent work of~\protect{\cite{charous2022stable}}.

It was originally studied around the same time in 1937 in two independent, pioneering works. 
Fisher~\protect{\cite{Fisher:1937}} studied a deterministic version of a stochastic model for the spread of a favored gene in a population in a one-dimensional habitat, with a ``logistic'' reaction term.
Kolmogorov, Petrowsky, and Piskunov provided a rigorous study of the two-dimensional equation and obtained some fundamental analytical results, with a general reaction term. We refer the reader to~\protect{\cite{Kolmogorov:1996}} for an English translation of their original work.

The Fisher--KPP equation can be used to model several phenomena in physics, chemistry, and biology. For instance, it can be used to describe biological population or chemical reaction dynamics with diffusion. It has also been used in the theory of combustion to study flame propagation and nuclear reactors. See~\protect{\cite[\S 13.2]{Murray:2002}} for a comprehensive review. 

\subsection{Boundary and initial conditions}

Here, we adopt the same boundary and initial conditions as in~\protect{\cite{charous2022stable}}. So the reaction rate is modeled as a random variable that follows a uniform law $ r \sim \cU \left[ 1/4, 1/2 \right] $. 
We consider the spatial domain: $ x \in [0, 40] $, and the time domain: $ t \in [0, 10] $. We impose homogeneous Neumann boundary conditions, i.e.,
\[
   \forall t \in [0, 10], \quad \frac{\partial w}{\partial x}(0,t) = 0, \quad \frac{\partial w}{\partial x}(40,t) = 0.
\]
These boundary conditions represent the physical condition of zero diffusive fluxes at the two boundaries.
The initial condition is ``stochastic'', of the form
\[
   w(x,0; \omega) = a(\omega) \, e^{-b(\omega) \, x^{2}},
\]
where $ a \sim \cU\left[1/5, \ 2/5\right] $ and $ b \sim \cU\left[1/10, \ 11/10\right] $.
The random variables $a$, $b$, and $r$ are all independent, and we consider $ N_{r} = 1000 $ realizations.

\subsection{Reference solution with the IMEX-CNLF method}

To obtain a reference solution, we use the implicit-explicit Crank--Nicolson leapfrog scheme (IMEX-CNLF) for time integration~\protect{\cite[Example~IV.4.3]{HundsdorferVerwer:2003}}. This scheme treats the linear diffusion term with Crank--Nicolson, an implicit method. In contrast, the nonlinear reaction term is treated explicitly with leapfrog, a numerical scheme based on the implicit midpoint method.

For the space discretization, we consider 1000 grid points in $ x $, while for the time discretization, we use $ 1601 $ points in time, so that the time step is $ \deltatime = 10/(1601 - 1) = 0.00625 $.

Let $ w^{(i)} $ denote the spatial discretization of the $ i $th realization.
At a given time $t$, each realization is stored as a column of our solution matrix, i.e.,
\[
   W(t) =
  \left[
  \begin{array}{cccc}
    \vertbar & \vertbar &   & \vertbar \\
       w^{(1)} &    w^{(2)} &  \cdots  &   w^{(N_{r})}  \\
    \vertbar & \vertbar &   & \vertbar 
  \end{array}
  \right].
\]
Moreover, let $ R_{\omega} $ be a diagonal matrix whose diagonal entries are the $ r_{(\omega)}^{(i)} $ coefficients for every realization indexed by $ i $, $ i = 1, 2, \ldots, N_{r} $. Indeed,
\begin{scriptsize}
\[
  \left[
  \begin{array}{cccc}
    \vertbar & \vertbar &   & \vertbar \\
       r_{(\omega)}^{(1)} w^{(1)} &    r_{(\omega)}^{(2)} w^{(2)} &  \cdots  &   r_{(\omega)}^{(N_{r})} w^{(N_{r})}  \\
    \vertbar & \vertbar &   & \vertbar 
  \end{array}
  \right]
  =
  \left[
  \begin{array}{cccc}
    \vertbar & \vertbar &   & \vertbar \\
       w^{(1)} &    w^{(2)} &  \cdots  &   w^{(N_{r})}  \\
    \vertbar & \vertbar &   & \vertbar 
  \end{array}
  \right] \cdot 
  \begin{bmatrix}
     r_{(\omega)}^{(1)} & & & \\
     & r_{(\omega)}^{(2)} & & \\
     & & \ddots & \\
     & & & r_{(\omega)}^{(N_{r})}
  \end{bmatrix} = W R_{\omega}.
\]
\end{scriptsize}

\begin{figure}[htbp]
    \centering
    \begin{minipage}[b]{0.33\textwidth}
        \centering
        \includegraphics[width=\textwidth]{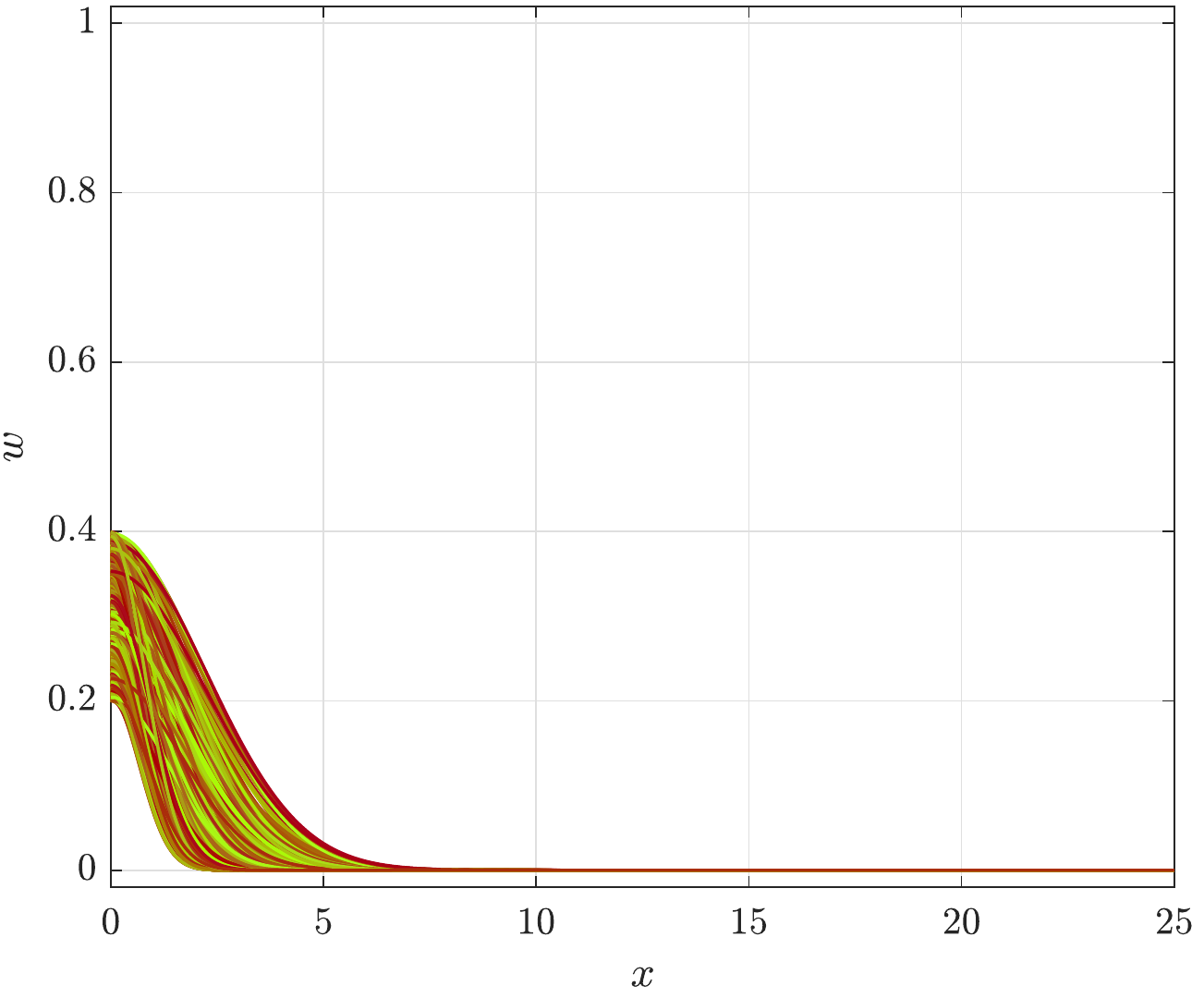}
        {\scriptsize (a)}
    \end{minipage}\hfill
    \begin{minipage}[b]{0.33\textwidth}
        \centering
        \includegraphics[width=\textwidth]{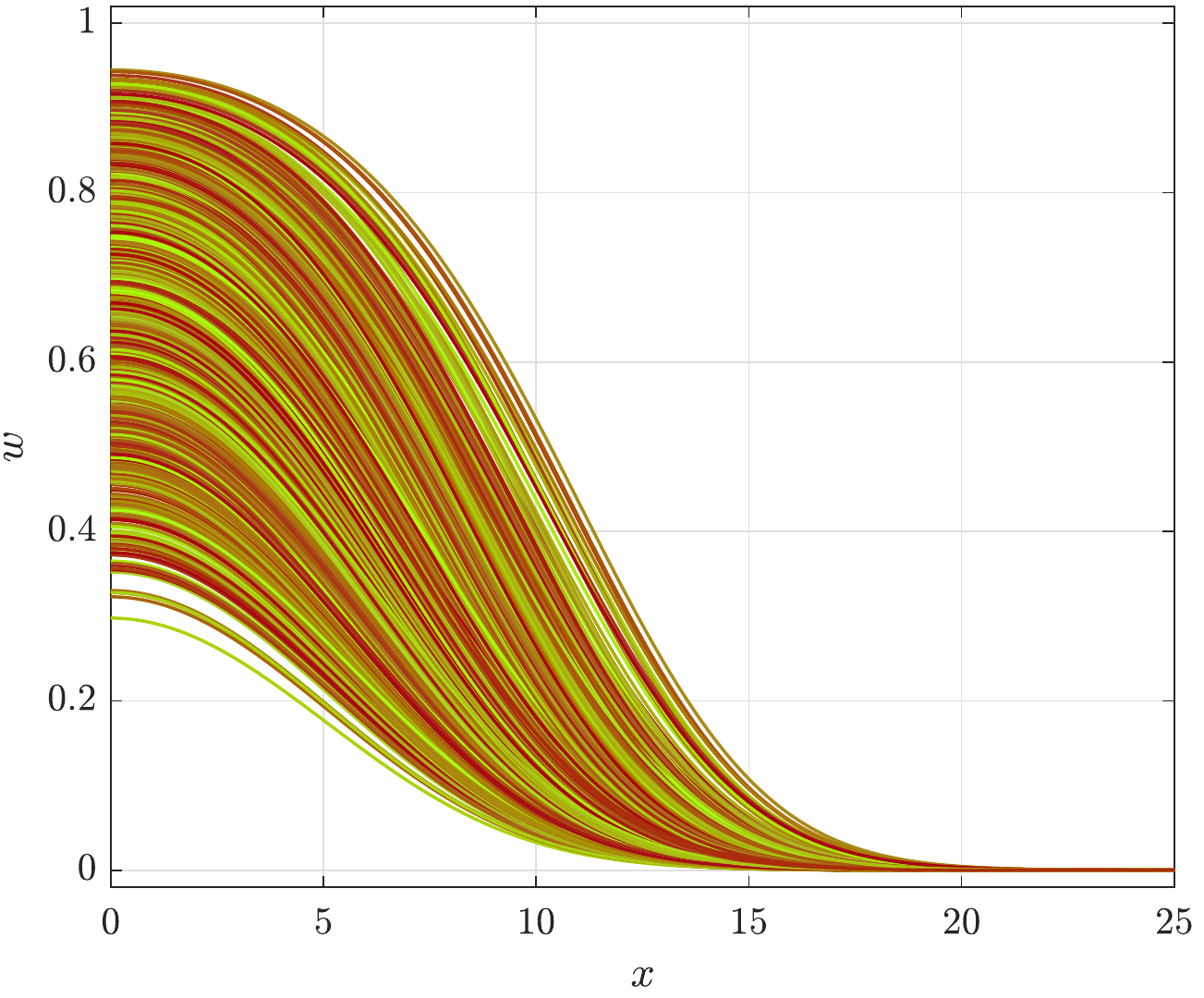}
        {\scriptsize (b)}
    \end{minipage}
    \begin{minipage}[b]{0.33\textwidth}
        \centering
        \includegraphics[width=\textwidth]{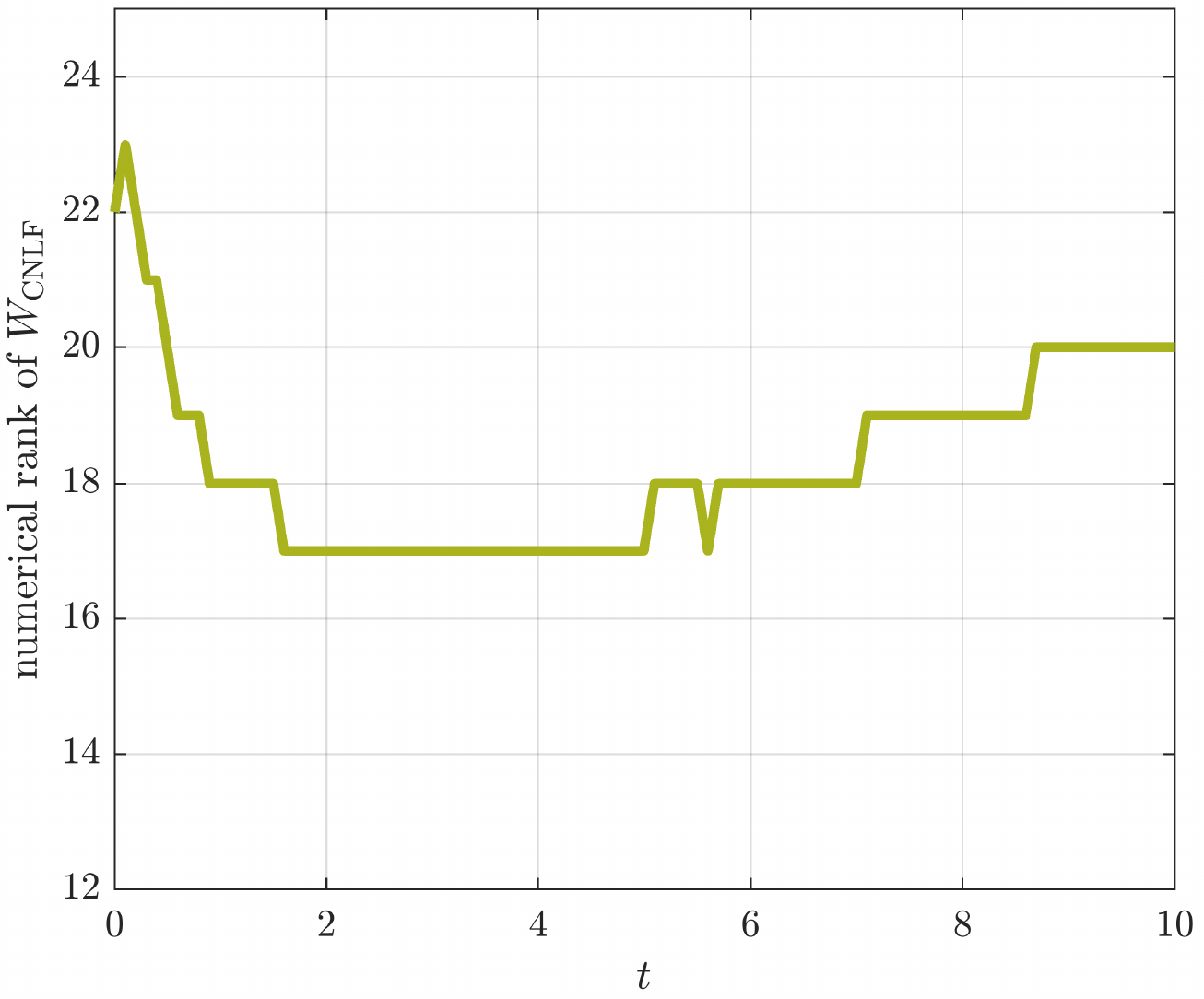}
        {\scriptsize (c)}
    \end{minipage}    
    \caption{Fisher--KPP reference solution computed with an IMEX-CNLF scheme. Panel (a): all the 1000 realizations at $ t = 0 $. Panel (b): all the 1000 realizations at $ t = 10 $. Panel (c): numerical rank history.}\label{fig:121314}
\end{figure}

The IMEX-CNLF scheme applied to~\eqref{eq:FKPP_equation} gives the algebraic equation 
\begin{equation}\label{eq:CNLS_alg_eq}
  (I - \deltatime A) W^{(n+1)} = (I + \deltatime A) W^{(n-1)} + 2 \deltatime \, W^{(n)} R_{\omega} - 2 \deltatime \, ( W^{(n)} )^{\circ 2} R_{\omega},
\end{equation}
where $ A $ is the matrix that discretizes the Laplacian with a second-order centered finite difference stencil and homogeneous Neumann boundary conditions, i.e.,
\begin{equation}\label{eq:discretized_laplacian_FKPP}
   A = \frac{1}{h_{x}^{2}}
   \begin{bmatrix}
      -2  &       2  &          &         &      \\
       1  &      -2  &       1  &         &      \\
          &  \ddots  &  \ddots  & \ddots  &      \\
          &          &       1  &     -2  &    1 \\
          &          &          &      2  &   -2
   \end{bmatrix}.
\end{equation}
For ease of notation, we call $ M_{\mathrm{m}} = I - \deltatime A $ and $ M_{\mathrm{p}} = I + \deltatime A $, so that~\eqref{eq:CNLS_alg_eq} becomes
\begin{equation}\label{eq:CNLS_alg_eq_2}
  M_{\mathrm{m}} W^{(n+1)} = M_{\mathrm{p}} W^{(n-1)} + 2 \deltatime \, W^{(n)} R_{\omega} - 2 \deltatime \, ( W^{(n)} )^{\circ 2} R_{\omega},
\end{equation}
Panels (a) and (b) of Figure~\ref{fig:121314} show the 1000 realizations at $ t = 0 $ and at $ t = 10 $, respectively.
Panel (c) reports on the numerical rank history. For computing the numerical rank, we use MATLAB's default tolerance, which in this case is about $ 10^{-11} $.


\subsection{LR-CNLF scheme}
To obtain a low-rank solver for the Fisher--KPP PDE, we proceed similarly as for the Allen--Cahn equation low-rank solution. We build a cost function $ F(W) $, so that its minimization gives the solution to~\eqref{eq:CNLS_alg_eq_2}.
\[
  F^{(n+1)}(W) = \frac{1}{2} \left\lVert M_{\mathrm{m}} W - M_{\mathrm{p}} W^{(n-1)} + 2 \deltatime \left(\big( W^{(n)} \big)^{\circ 2} -  W^{(n)}  \right) R_{\omega} \right\rVert_{\F}^{2}.
\]
Developing and keeping only the terms that depend on $ W $, we get the cost function:
\begin{small}
\begin{equation}\label{eq:FKPP_cost_function}
  F(W) = \frac{1}{2} \trace\!\big( W\tr \! M_{\mathrm{m}}\tr \! M_{\mathrm{m}} W \big) - \trace\!\big((W^{(n-1)})\tr \! M_{\mathrm{p}}\tr \! M_{\mathrm{m}} W\big) + 2 \deltatime \trace\!\left(\left(\big( W^{(n)} \big)^{\circ 2} -  W^{(n)}  \right)\tr \! M_{\mathrm{m}} W R_{\omega}\right).
\end{equation}
\end{small}
More details are given in Appendix~\ref{app:discretization_FKPP}.

\subsection{Numerical experiments}

We monitor the following quantities:
\begin{itemize}
\item the numerical rank of $ W_{\mathrm{RTR}} $ (the numerical rank of $ W_{\mathrm{CNLF}} $ is also plotted as reference);
\item the $ L^{2} $-norm of the error $ \| w_{\mathrm{RTR}} - w_{\mathrm{CNLF}} \|_{L^{2}(\Omega)} = \| W_{\mathrm{RTR}} - W_{\mathrm{CNLF}} \|_{\F} \cdot \sqrt{h_{x}} $.
\end{itemize}

\begin{figure}[htbp]
    \centering
    \begin{minipage}[b]{0.485\textwidth}
        \centering
        \includegraphics[width=\textwidth]{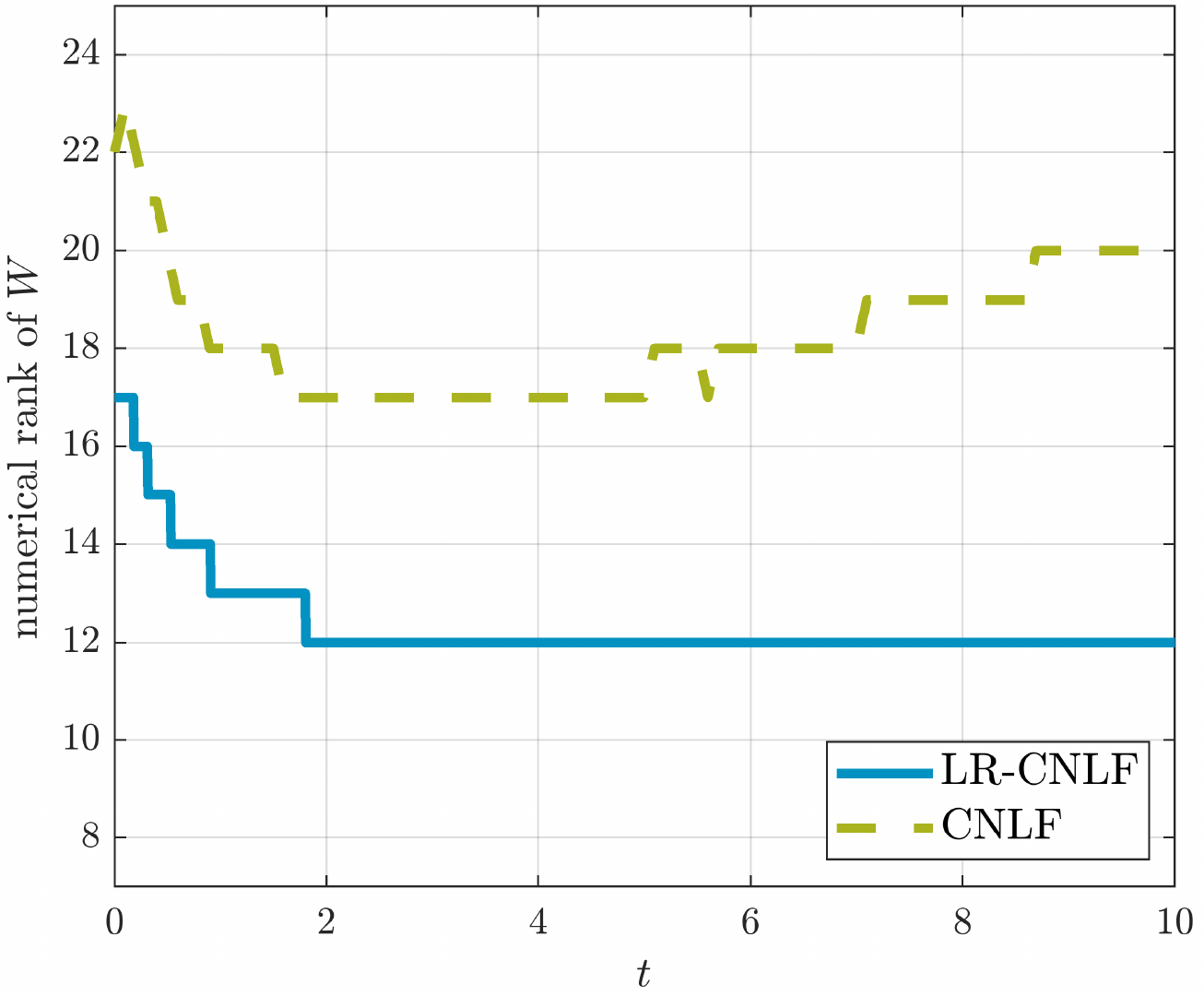}
        {\scriptsize (a)}
    \end{minipage}\hfill
    \begin{minipage}[b]{0.500\textwidth}
        \centering
        \includegraphics[width=\textwidth]{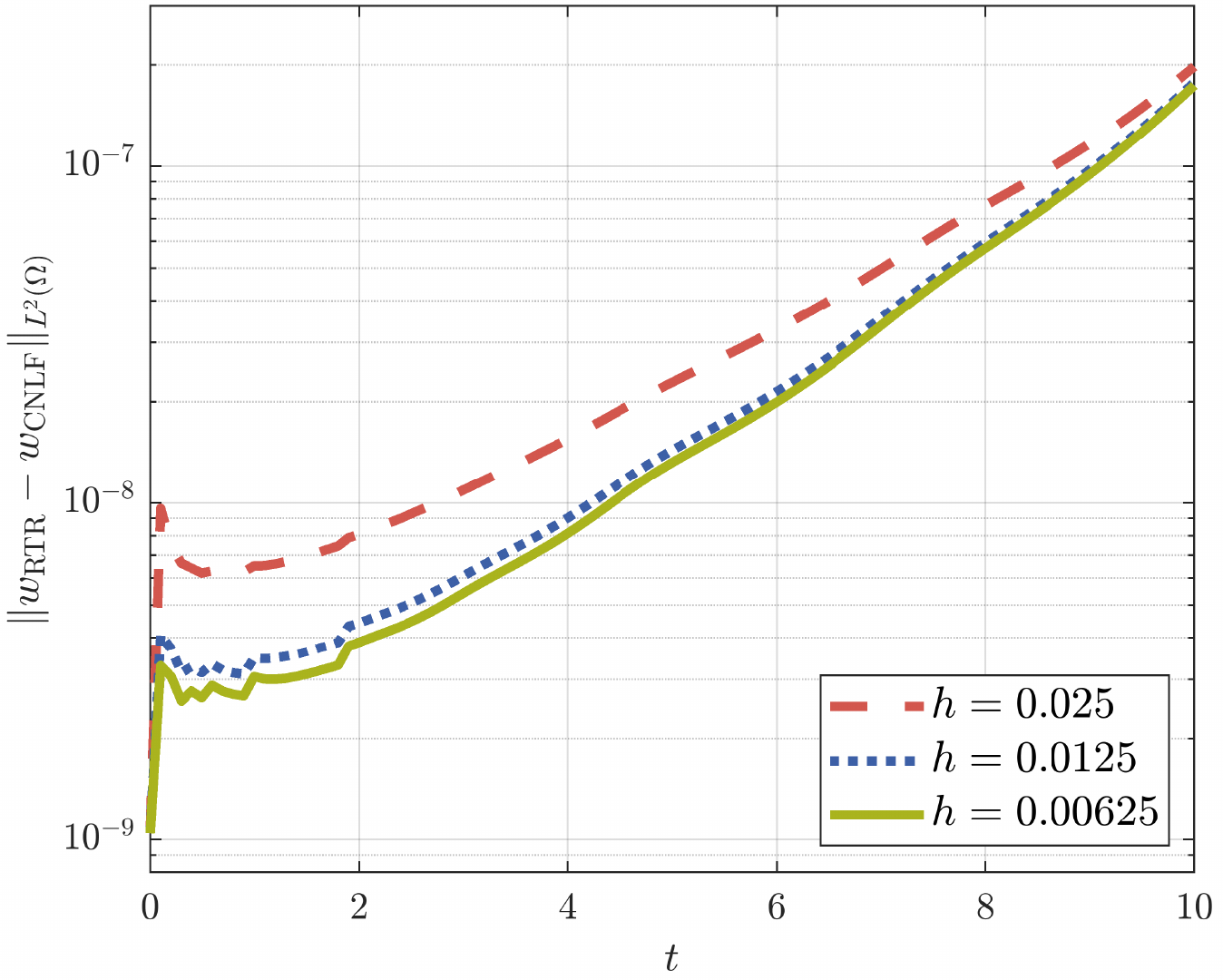}
        {\scriptsize (b)}
    \end{minipage}
    \caption{Panel (a): rank history for the LR-CNLF method compared to the reference solution, for $ \deltatime = 0.00625 $. Panel (b): discrete $ L^{2} $-norm of the error versus time, for several $ \deltatime $.}\label{fig:1516}
\end{figure}

As it was done in the previous section for the reference solution, here we also consider 1000 realizations.
We apply our technique with rank adaption, with tolerance for rank truncation of $ 10^{-8} $. The inner PrecRTR is halted once the gradient norm is less than $ 10^{-8} $. Figure~\ref{fig:1516} reports on the numerical experiments.

\section{Numerical experiments for LYAP and NPDE} \label{sec:numerical_experiments}

This section focuses on the numerical properties of PrecRTR, our preconditioned RTR on the manifold of fixed-rank matrices on the variational problems from~\cite{Sutti_V:2021}, and that were recalled in Sect.~\ref{sec:problem_setting}. These are large-scale finite-dimensional optimization problems arising from the discretization of infinite-dimensional problems. These problems have been used as benchmarks in nonlinear multilevel algorithms, as seen in~\cite{Henson:2003,Gratton:2008,Wen:2009}. For further information on the theoretical aspects of variational problems, we recommend consulting~\cite{Brenner:2007,LeDret:2016}.

We consider two scenarios: in the first one, we let Manopt automatically take care of $ \bar{\Delta} $. In the second one, we fix $ \bar{\Delta} = 0.5 $.
The tolerance on the norm of the gradient in the trust-region method is set to $ 10^{-12} $, and we set the maximum number of outer iterations $ n_{\text{max outer}} = 300 $.

We report on the behavior of PrecRTR by monitoring the Frobenius norm of the normalized Riemannian gradient:
\[
   \textrm{R-grad}(i) \coloneqq \| \xi^{(i)} \|_{\F}/\| \xi^{(0)} \|_{\F}.
\]

\subsection{Tables}

Tables~\ref{tab:comparison_LYAP} and~\ref{tab:comparison_NPDE} report the convergence results for PrecRTR. 
CPU time is in seconds, and for each line in the tables, the experiments are repeated and averaged over 10 times.

To assess the accuracy of the solutions obtained for the Lyapunov equation, we also use the standard residual
\[
   r(W ) \coloneqq \| AW  + WA - \Gamma  \|_{\F}.
\]

When the maximum number of outer iterations $ n_{\text{max outer}} = 300 $ is reached, we indicate this in bold text.
We also set a limit on the quantity $ \sum n_{\mathrm{inner}} $: the inner tCG solver is stopped when $ \sum n_{\mathrm{inner}} $ first exceeds $ 30\,000 $. This is highlighted by the bold text in the tables below.

\begin{table}[htbp]
   \caption{Preconditioned RTR for the LYAP problem.}
   \label{tab:comparison_LYAP}
   \begin{center}
   \resizebox{\textwidth}{!}{%
      \begin{tabular}{cr|ccc|ccc}  
         \cmidrule[1pt](lr){3-8}
         \multicolumn{2}{c}{} &  \multicolumn{3}{c}{Rank 5} & \multicolumn{3}{c}{Rank 10}  \\
         \midrule
          $ \ell $  &   size     &  time    &  $ \| \xi^{(\mathrm{end})} \|_{F} $  &   $ r(W^{(\mathrm{end})}) $   &  time    &  $ \| \xi^{(\mathrm{end})} \|_{F} $  &   $ r(W^{(\mathrm{end})}) $   \\
         \midrule
            10  &   1\,048\,576 &        0.21  &  $ 1.0150 \times 10^{-13} $  &  $ 9.7480 \times 10^{-8} $  &        0.85  &  $ 6.2481 \times 10^{-14} $  &  $ 4.2204 \times 10^{-11} $  \\
            11  &   4\,194\,304 &        0.49  &  $ 2.9645 \times 10^{-14} $  &  $ 4.8741 \times 10^{-8} $  &        1.53  &  $ 5.7690 \times 10^{-13} $  &  $ 2.0374 \times 10^{-11} $  \\
            12  &  16\,777\,216 &        1.01  &  $ 3.8413 \times 10^{-14} $  &  $ 2.4371 \times 10^{-8} $  &        2.93  &  $ 1.0921 \times 10^{-13} $  &  $ 1.0478 \times 10^{-11} $  \\
            13  &  67\,108\,864 &        1.56  &  $ 7.3017 \times 10^{-14} $  &  $ 1.2185 \times 10^{-8} $  &        5.74  &  $ 1.3556 \times 10^{-13} $  &  $ 5.2396 \times 10^{-12} $  \\
            14  & 268\,435\,456 &        3.80  &  $ 1.5082 \times 10^{-13} $  &  $ 6.0927 \times 10^{-9} $  &        10.87  &  $ 9.3753 \times 10^{-14} $  &  $ 2.6045 \times 10^{-12} $  \\
            15  & 1\,073\,741\,824 &     7.48  &  $ 2.7525 \times 10^{-13} $  &  $ 3.0464 \times 10^{-9} $  &        25.02  &  $ 2.4835 \times 10^{-13} $  &  $ 1.3177 \times 10^{-12} $   \\
         \bottomrule
      \end{tabular}}
   \end{center}
\end{table}

\begin{table}[htbp]
   \caption{Effect of preconditioning: dependence on $ \ell $ for LYAP.}
   \label{tab:size_LYAP}
      \begin{center}
      \resizebox{\textwidth}{!}{
      \begin{tabular}{c|c|cccccc|cccccc}  
         \cmidrule[1pt](lr){3-14}
            \multicolumn{2}{c}{} &  \multicolumn{6}{c}{Rank 5} &  \multicolumn{6}{c}{Rank 10}  \\
         \midrule
            Prec. &  $ \ell $   &    10  &     11  &     12  &     13  &   14  &   15 &     10  &     11  &     12  &     13  &   14  &   15 \\
         \midrule
         \midrule
         \multirow{3}{*}{No}
             &  $ n_{\mathrm{outer}} $       &    51  &     54  &     61  &           59  &         162   &           92  &  \textbf{300}  &         103  &           61 &                 63           &           62   &    59 \\
             &  $ \sum n_{\mathrm{inner}} $  &  4561  &   9431  &  21066  &  \textbf{36556}  &  \textbf{30069}  &  \textbf{30096}  &      27867  & \textbf{30025}  &  \textbf{33818} &       \textbf{45760}  &   \textbf{44467}  &      \textbf{38392}  \\
             &  $ \max n_{\mathrm{inner}} $  &  1801  &   3191  &   7055  &  9404  &  1194  &   1851   &   2974  &       3385  &       8894  &        24367  &        24537  &   25013   \\
         \midrule
         \multirow{3}{*}{Yes}
             &  $ n_{\mathrm{outer}} $       &    41  &     45  &    50  &    52  &    56  &    60   &    44  &    64  &    62  &   53  &   56  &   56  \\
             &  $ \sum n_{\mathrm{inner}} $  &    44  &     45  &    50  &    52  &    56  &    60   &    69  &   104  &    82  &   60  &   69  &   56  \\
             &  $ \max n_{\mathrm{inner}} $  &     4  &      1  &     1  &     1  &     1  &     1   &     9  &     9  &     8  &    8  &    8  &    1  \\
         \bottomrule
      \end{tabular}}
   \end{center}
\end{table}

Tables~\ref{tab:size_LYAP} and~\ref{tab:size_NPDE} report on the effect of preconditioning as $ \ell $ increases for LYAP and NPDE, respectively. 
The reductions in the number of iterations of the inner tCG between the non-preconditioned and the preconditioned versions are impressive.
Moreover, for the preconditioned method, both tables demonstrate that $ n_{\mathrm{outer}} $ and $ \sum n_{\mathrm{inner}} $ depend (quite mildly) on $ \ell $, while $ \max n_{\mathrm{inner}} $ is basically constant.

For NPDE, in both the non-preconditioned and preconditioned methods, the numbers of iterations are typically higher than those for the LYAP problem, which is plausibly due to the nonlinearity of the problem.

Tables~\ref{tab:rank_fixed_level_12_LYAP} and~\ref{tab:rank_fixed_level_12_NPDE} reports the results for varying rank $ r $ and fixed problem size $ \ell = 12 $. The stopping criteria are the same as above. It is remarkable that, for PrecRTR for the LYAP problem, all three monitored quantities basically do not depend on the rank $ r $. For NPDE, there is some more, but still quite mild, dependence on $ r $.

\begin{table}[htbp]
   \caption{Effect of preconditioning: dependence on $ r $ with fixed size $ \ell = 12 $, for LYAP.}
   \label{tab:rank_fixed_level_12_LYAP}
      \begin{center}
      \begin{tabular}{c|c|cccccc}  
         \cmidrule[1pt](lr){3-8}
                \multicolumn{2}{c}{}         &   \multicolumn{6}{c}{Rank}   \\
         \midrule
            Prec. &       iterations         &     1   &     2   &       5  &          10  &        15   &        20  \\
         \midrule
         \midrule
         \multirow{3}{*}{No}
             &  $ n_{\mathrm{outer}} $       &     53  &     53  &     61  &           61  &  \textbf{300}  &          62  \\
             &  $ \sum n_{\mathrm{inner}} $  &  17650  &  18775  &  21066  &  \textbf{33818}  &      12816  &  \textbf{33292} \\
             &  $ \max n_{\mathrm{inner}} $  &   6276  &   7225  &   7055  &         8894  &       3794  &        6928  \\
         \midrule
         \multirow{3}{*}{Yes}
             &  $ n_{\mathrm{outer}} $       &     51  &     51  &    50  &            49  &         49  &       48  \\
             &  $ \sum n_{\mathrm{inner}} $  &     51  &     51  &    50  &            49  &         49  &       48  \\
             &  $ \max n_{\mathrm{inner}} $  &      1  &      1  &     1  &             1  &          1  &        1  \\
         \bottomrule
      \end{tabular}
   \end{center}
\end{table}

\begin{table}[htbp]
   \caption{Preconditioned RTR for the NPDE problem.}
   \label{tab:comparison_NPDE}
   \begin{center}
   \resizebox{\textwidth}{!}{%
      \begin{tabular}{cr|ccc|ccc}  
         \cmidrule[1pt](lr){3-8}
         \multicolumn{2}{c}{} &  \multicolumn{3}{c}{Rank 5} & \multicolumn{3}{c}{Rank 10}  \\
         \midrule
          $ \ell $  &   size     &  time    &  $ \| \xi^{(\mathrm{end})} \|_{F} $  &   $ r(W^{(\mathrm{end})}) $   &  time    &  $ \| \xi^{(\mathrm{end})} \|_{F} $  &   $ r(W^{(\mathrm{end})}) $   \\
         \midrule \midrule
         \multicolumn{8}{c}{Rank 5}  \\
         \midrule
            10  &   1\,048\,576 &        0.45  &  $ 2.0719 \times 10^{-14} $  &  $ 1.5614 \times 10^{-5} $  &        1.17  &  $ 1.7303 \times 10^{-14} $  &  $ 1.8660 \times 10^{-7} $  \\
            11  &   4\,194\,304 &        0.89  &  $ 2.7106 \times 10^{-14} $  &  $ 7.8072 \times 10^{-6} $  &        2.10  &  $ 6.0181 \times 10^{-14} $  &  $ 9.3301 \times 10^{-8} $  \\
            12  &  16\,777\,216 &        1.65  &  $ 5.2974 \times 10^{-14} $  &  $ 3.9036 \times 10^{-6} $  &        4.73  &  $ 5.9537 \times 10^{-14} $  &  $ 4.6650 \times 10^{-8} $  \\
            13  &  67\,108\,864 &        2.84  &  $ 1.2492 \times 10^{-13} $  &  $ 1.9518 \times 10^{-6} $  &        8.91  &  $ 1.1536 \times 10^{-13} $  &  $ 2.3325 \times 10^{-8} $  \\
            14  & 268\,435\,456 &        5.89  &  $ 2.4349 \times 10^{-13} $  &  $ 9.7591 \times 10^{-7} $  &       19.67  &  $ 2.6992 \times 10^{-13} $  &  $ 1.1663 \times 10^{-8} $  \\
            15  & 1\,073\,741\,824 &    12.96  &  $ 6.4490 \times 10^{-13} $  &  $ 4.8796 \times 10^{-7} $  &       45.71  &  $ 5.8336 \times 10^{-13} $  &  $ 5.8313 \times 10^{-9} $  \\
         \bottomrule
      \end{tabular}}
   \end{center}
\end{table}

\begin{table}[htbp]
   \caption{Effect of preconditioning: dependence on $ \ell $ for NPDE.}
   \label{tab:size_NPDE}
      \begin{center}
      \resizebox{\textwidth}{!}{
      \begin{tabular}{c|c|cccccc|cccccc}  
         \cmidrule[1pt](lr){3-14}
            \multicolumn{2}{c}{} &  \multicolumn{6}{c}{Rank 5} &  \multicolumn{6}{c}{Rank 10}  \\
         \midrule
            Prec. &  $ \ell $   &    10  &     11  &     12  &    13  &   14  &   15 &     10 &    11  &    12  &     13  &   14  &   15 \\
         \midrule
         \midrule
         \multirow{3}{*}{No}
             &  $ n_{\mathrm{outer}} $       &    53  &     57  &     61  &     79 &   68  &  68  &  63 &  87 & 76 &  68  &   62   &    65    \\
             &  $ \sum n_{\mathrm{inner}} $  &  4603  &   9505  &  13817  &  \textbf{41144} &  \textbf{47186}  &  \textbf{38079}  &      6610  &  \textbf{38858}  &    \textbf{30567}  &  \textbf{31028}  & \textbf{39803} & \textbf{39337} \\
             &  $ \max n_{\mathrm{inner}} $  &  2022  &   3595  &   7735  &  14195 &   28410  &  32433  &   1487  &   11550  &   6035  &  10598  &    22468   &   30118  \\
         \midrule
         \multirow{3}{*}{Yes}
             &  $ n_{\mathrm{outer}} $       &    50  &     56  &     61  &    63  &   65  &     66  &        53  &       58  &       63  &     69  &   69  &  71 \\
             &  $ \sum n_{\mathrm{inner}} $  &    57  &     64  &     69  &    72  &   74  &     75  &        78  &       84  &       90  &     98  &   97  &  100  \\
             &  $ \max n_{\mathrm{inner}} $  &     6  &      7  &      7  &     7  &    7  &      7  &        10  &       10  &       11  &     11  &   10  &   10  \\
         \bottomrule
      \end{tabular}}
   \end{center}
\end{table}

\begin{table}[htbp]
   \caption{Effect of preconditioning: dependence on $ r $ for NPDE with fixed level $ \ell = 12 $.}
   \label{tab:rank_fixed_level_12_NPDE}
      \begin{center}
      \begin{tabular}{c|c|cccccc}  
         \cmidrule[1pt](lr){3-8}
                \multicolumn{2}{c}{}         &   \multicolumn{6}{c}{Rank}   \\
         \midrule
            Prec. &       iterations         &     1  &      2  &       5  &       10  &       15  &           20  \\
         \midrule
         \midrule
         \multirow{3}{*}{No}
             &  $ n_{\mathrm{outer}} $       &    59  &     57  &       61  &          76  &           62  &           60  \\
             &  $ \sum n_{\mathrm{inner}} $  &  9183  &  16044  &    13826  &  \textbf{30567} &  \textbf{61339}  &  \textbf{31192}  \\
             &  $ \max n_{\mathrm{inner}} $  &  3569  &   4642  &     7744  &         6035 &        31627  &         8540  \\
         \midrule
         \multirow{3}{*}{Yes}
             &  $ n_{\mathrm{outer}} $       &    59  &    61  &       61  &         63 &       60  &          61  \\
             &  $ \sum n_{\mathrm{inner}} $  &    78  &    90  &       69  &         90 &       90  &         104  \\
             &  $ \max n_{\mathrm{inner}} $  &    11  &    10  &        7  &         11 &       11  &          13  \\
         \bottomrule
      \end{tabular}
   \end{center}
\end{table}

\section{Conclusions and outlook} \label{sec:conclusions}

In this paper, we have shown how to combine an efficient preconditioner with optimization on low-rank manifolds. Unlike classical Lyapunov solvers, our optimization strategy can treat nonlinearities. Moreover, compared to iterative methods that perform rank-truncation at every step, our approach allows for much larger time steps as it does not need to satisfy a fixed-point Lipschitz restriction. We illustrated the application of this technique to two time-dependent nonlinear PDEs --- the Allen--Cahn and the Fisher--KPP equations. In addition, the numerical experiments for two variational problems demonstrate the efficiency in computing good low-rank approximations with a number of tCG iterations in the trust region subsolver which is almost independent of the problem size.
Future research may focus on higher-order methods, such as more accurate implicit methods. Additionally, we may explore higher-dimensional problems, problems in biology, and stochastic PDEs.

\section*{Acknowledgments}
The work of the first author was supported by the National Center for Theoretical Sciences and the Ministry of Science and Technology of Taiwan (R.O.C.), under the contracts 111-2124-M-002-014 and 112-2124-M-002-009.

\appendix

\section{Low-rank formats for the Allen--Cahn equation}\label{app:discretization_ACE}
\subsection{Objective functional}\label{app:discr_objective}

Discretizing~\eqref{eq:ACE_var_pb} similarly as in~\protect{\cite[\S 5.2.1]{Sutti_V:2021}}, we obtain
\begin{equation}\label{eq:discr_objective}
   F = h_{x}^{2} \sum_{i,j=0}^{2^{\ell} - 1} \left( \frac{\varepsilon \deltatime}{2} ( \partial w_{x_{ij}}^{2} + \partial w_{y_{ij}}^{2} ) + \frac{1-\deltatime}{2} \, w_{ij}^{2} + \frac{\deltatime}{4} \, w_{ij}^{4} - \widetilde{w}_{ij} w_{ij} \right).
\end{equation}
To obtain the factored format of the discretized objective functional, we consider the factorizations $ W = U\Sigma V\tr $, and $ \widetilde{W} = \widetilde{U} \widetilde{\Sigma} \widetilde{V}\tr $.

The first term and the fourth term in~\eqref{eq:discr_objective} have the same factorized form as those seen in~\protect{\cite[\S 5.2.1]{Sutti_V:2021}}.
The only slight change is due to the periodic boundary conditions adopted here. As a consequence, the matrix $ L $ that discretizes the first-order derivatives\footnote{Sometimes known as \emph{forward difference matrix}.} with periodic boundary conditions becomes
\begin{equation*}
   L = \frac{1}{h_{x}}
   \begin{bmatrix}
       -1 &   1 &          &         &    \\
          &  -1 &        1 &         &    \\ 
          &     &   \ddots &  \ddots &    \\
        1 &     &          &      -1 &  1
   \end{bmatrix}.
\end{equation*}
Note the presence of the unitary coefficient in the lower-left corner. The reader can easily verify that $ A = L\tr L $, where $ A $ is the matrix~\eqref{eq:discretized_laplacian}.
We recall that, given this matrix, the first-order derivatives of $ W $ can be computed as
\[
   \partial W_{x} = L W  \quad \text{and} \quad  \partial W_{y} = W L\tr.
\]
For the second term in~\eqref{eq:discr_objective}, we have
\[
   \sum_{i,j} \frac{1-\deltatime}{2} \, w_{ij}^{2} = \frac{1-\deltatime}{2} \, \trace(W\tr W) = \frac{1-\deltatime}{2} \, \| \Sigma \|^{2}_{\F},
\]
For the third term, it is easier to consider the full-rank format
\begin{equation}\label{eq:modified_third_term}
    \frac{\deltatime}{4} \sum_{i,j} w_{ij}^{4}.
\end{equation}
We call $ \widetilde{G} = \widetilde{U}\widetilde{\Sigma} $ and $ G = U \Sigma $.
Finally, the discretized objective functional in factorized matrix form is
\[
    F = h_{x}^{2} \left( \frac{\varepsilon \deltatime}{2} \left( \| (LU) \Sigma \|^{2}_{\F} + \| (LV) \Sigma \|^{2}_{\F} \right) + \frac{1-\deltatime}{2} \| \Sigma \|^{2}_{\F} + \frac{\deltatime}{4} \sum_{i,j} w_{ij}^{4} - \trace\!\big( (\widetilde{G}\tr G) (V\tr \widetilde{V}) \big) \right).
\]

Table~\ref{tab:complexities_ACE_cost} summarizes the asymptotic complexities for the ACE cost function. In this and the following similar tables, we indicate the sizes of the matrices in the order in which they appear in the product. If all the matrices in a term are the same size, we only indicate that once. Matrices without any specific structure are stored as dense, unless otherwise specified.

\begin{table}[htbp]
   \caption{Asymptotic complexities for ACE cost function.}
   \label{tab:complexities_ACE_cost}
   \begin{center}
      \begin{tabular}{c|c|c|c}  
         \toprule
            Product  & Factor sizes &  Notes on structure and storage  &  Cost   \\
         \midrule
         \midrule
            $ \| (LU) \Sigma \|^{2}_{\F} $  &  $ n \times n $,  $ n \times r $,  $ r \times r $  &  $ L $ sparse banded, $ \Sigma $ sparse diagonal  &  $ \cO(nr) $ \\
            $ \| (LV) \Sigma \|^{2}_{\F} $  &  $ n \times n $,  $ n \times r $,  $ r \times r $  &  $ L $ sparse banded, $ \Sigma $ sparse diagonal  &  $ \cO(nr) $ \\
            $ \| \Sigma \|^{2}_{\F} $  &  $ r \times r $  &  $ \Sigma $ sparse diagonal  &  $ \cO(r) $ \\
            $ \sum_{i,j} w_{ij}^{4} $  &  $ n \times n $  &    &  $ \cO(n^{2}) $ \\
          \midrule
          \multicolumn{3}{l}{$ \trace\!\big( (\widetilde{G}\tr\! G) (V\tr \widetilde{V}) \big) $} &  $ \cO(nr^{2} + r^{3}) $  \\
          \midrule
            $ \widetilde{G}\tr\! G $  & $ n \times r $,  $ n \times r $ &  &  $ \cO(nr^{2}) $ \\
            $ V\tr \widetilde{V} $  &  $ n \times r $,  $ n \times r $ &  &  $ \cO(nr^{2}) $ \\
            $ (\widetilde{G}\tr\! G) (V\tr \widetilde{V}) $  & $ r \times r $, $ r \times r $  &  &  $ \cO(r^{3}) $ \\
            $ \trace(\cdot) $  & $ r \times r $  &  &  $ \cO(r) $ \\
         \bottomrule
      \end{tabular}
   \end{center}
\end{table}

\subsection{Gradient}

The gradient of $ \cF $~\eqref{eq:ACE_var_pb} is the variational derivative
\[
   \frac{\delta \cF}{\delta w} = -\varepsilon \deltatime \Delta w + (1 - \deltatime) \, w + \deltatime \, w^{3} - \widetilde{w}.
\]
The discretized Euclidean gradient in matrix form is given by
\[
    G = h_{x}^{2} \left( - \varepsilon \deltatime ( A W + W A ) + (1-\deltatime) W + \deltatime W^{\circ 3} - \widetilde{W} \right),
\]
with $ A $ as in~\eqref{eq:discretized_laplacian}.

For the term $ W^{\circ 2} = W \odot W $, we perform the element-wise multiplication in factorized form as explained in \protect{\cite[§7]{Kressner:2014}} and store the result in the format $ U_{\circ 2} \Sigma_{\circ 2} V_{\circ 2}\tr $, i.e.,
\[
   W^{\circ 2} = W \odot W = ( U \kt U ) ( \Sigma \otimes \Sigma ) (V \kt V)\tr = U_{\circ 2} \Sigma_{\circ 2} V_{\circ 2}\tr,
\]
where $ \kt $ denotes a transposed variant of the Khatri--Rao product.
Then for $ W^{\circ 3} = W \odot W \odot W $ we consider the factorized format
\[
   W^{\circ 3} = W^{\circ 2} \odot W = ( U_{\circ 2} \kt U ) ( \Sigma_{\circ 2} \otimes \Sigma ) (V_{\circ 2} \kt V)\tr = U_{\circ 3}\Sigma_{\circ 3}V_{\circ 3}\tr,
\]

Substituting the formats $ W = U\Sigma V\tr $, $ W^{\circ 3} = U_{\circ 3}\Sigma_{\circ 3}V_{\circ 3}\tr $, and $ \widetilde{W} = \widetilde{U} \widetilde{\Sigma} \widetilde{V}\tr $, we get the factorized form of the Euclidean gradient $ G = U_{G} \Sigma_{G} V_{G}\tr $, where
\[
   U_{G} = \left[ \left( -\varepsilon \deltatime A + (1-\deltatime) I \right) U \quad U \quad U_{\circ 3} \quad \widetilde{U} \right],
\]
\[
   \Sigma_{G} = h_{x}^{2} \blkdiag \! \left( \Sigma, \ (-\varepsilon \deltatime) \Sigma, \ \deltatime \Sigma_{\circ 3}, \ -\widetilde{\Sigma} \right),
\]
and
\[
   V_{G} = \left[ V \quad A V \quad V_{\circ 3} \quad \widetilde{V} \right].
\]
The gradient $ G $ is an augmented matrix, analogously to the discretized gradient in factored format for the ``NPDE'' problem (see~\protect{\cite[\S 5.2.2]{Sutti_V:2021}}). The operations needed to form $ G $ are summarized in Table~\ref{tab:complexities_ACE_gradient}.

\begin{table}[htbp]
   \caption{Asymptotic complexities for ACE gradient.}
   \label{tab:complexities_ACE_gradient}
   \begin{center}
      \begin{tabular}{c|c|c|c}  
         \toprule
            Product  & Factor sizes &  Notes on structure and storage  &  Cost   \\
         \midrule
         \midrule
            $ \left( -\varepsilon \deltatime A + (1-\deltatime) I \right) U $  &  $ n \times n $, $ n \times r $  & $ A $, $ I $ sparse banded &  $ \cO(nr) $  \\
          \midrule
            $ A V $  &  $ n \times n $, $ n \times r $  & $ A $ sparse banded &  $ \cO(nr) $  \\
          \midrule
            $  U_{\circ 2} = U \kt U $  &  $ n \times r $  &  &  $ \cO(nr^{2}) $  \\
            $  \Sigma_{\circ 2} = \Sigma \otimes \Sigma $  &  $ r \times r $  &  $ \Sigma $ sparse diagonal  &  $ \cO(r^{2}) $  \\
            $  V_{\circ 2} = V \kt V $  &  $ n \times r $  &  &  $ \cO(nr^{2}) $  \\
          \midrule
            $  U_{\circ 3} = U_{\circ 2} \kt U $  &  $ n \times r^{2} $, $ n \times r $   &  &   $ \cO(nr^{3}) $  \\
            $  \Sigma_{\circ 3} = \Sigma_{\circ 2} \otimes \Sigma $  &  $ r^{2} \times r^{2} $, $ r \times r $  & $ \Sigma_{\circ 2} $, $ \Sigma $ sparse diagonal  &  $ \cO(r^{3}) $  \\
            $  V_{\circ 3} = V_{\circ 2} \kt V $  &  $ n \times r^{2} $, $ n \times r $ &  &  $ \cO(nr^{3}) $   \\
          \bottomrule
      \end{tabular}
   \end{center}
\end{table}

\subsection{Hessian}

The discretized Euclidean Hessian is (compare~\protect{\cite[\S 7.4.2.3]{Sutti:2020}})
\[
   H_{W}[\eta] = h_{x}^{2} \left( -\varepsilon \deltatime ( A \eta + \eta A ) + 3 \deltatime \, W^{\circ 2} \odot \eta + (1-\deltatime) \eta \right).
\]
The factored form of the discretized Euclidean Hessian is $ H_{W}[\eta] = U_{H_{W}[\eta]} S_{H_{W}[\eta]} V_{H_{W}[\eta]}\tr $, where
\[
   U_{H_{W}[\eta]} = \left[ \left( -\varepsilon \deltatime A + (1-\deltatime) I \right) U_{\eta} \quad U_{\eta} \quad U_{\odot} \right],
\]
\[
   S_{H_{W}[\eta]} = h_{x}^{2} \blkdiag\!\left( S_{\eta}, \ (-\varepsilon \deltatime) S_{\eta}, \ 3\deltatime \Sigma_{\odot} \right),
\]
\[
   V_{H_{W}[\eta]} = \left[ V_{\eta} \quad A V_{\eta} \quad V_{\odot} \right].
\]
where $ \eta = U_{\eta} S_{\eta} V_{\eta}\tr $ is a tangent vector in $ \mathrm{T}_{W}\cMr $, and $ W^{\circ 2} \odot \eta = U_{\odot} \Sigma_{\odot} V_{\odot}\tr $.

\begin{table}[htbp]
   \caption{Asymptotic complexities for ACE Hessian.}
   \label{tab:complexities_ACE_Hessian}
   \begin{center}
      \begin{tabular}{c|c|c|c}  
         \toprule
            Product  & Factor sizes &  Notes on structure and storage  &  Cost   \\
         \midrule
         \midrule
            $ \left( -\varepsilon \deltatime A + (1-\deltatime) I \right) U_{\eta} $  &  $ n \times n $, $ n \times r $  & $ A $ sparse banded, $ I $ sparse diagonal  &  $ \cO(nr) $  \\
          \midrule
            $ A V_{\eta} $  &  $ n \times n $, $ n \times r $  & $ A $ sparse banded &  $ \cO(nr) $  \\
          \midrule
            $  U_{\circ 2} = U \kt U $  &  $ n \times r $  &  &  $ \cO(nr^{2}) $  \\
            $  \Sigma_{\circ 2} = \Sigma \otimes \Sigma $  &  $ r \times r $  &  $ \Sigma $ sparse diagonal  &  $ \cO(r^{2}) $  \\
            $  V_{\circ 2} = V \kt V $  &  $ n \times r $  &  &  $ \cO(nr^{2}) $  \\
          \midrule
            $  U_{\odot} = U_{\circ 2} \kt U_{\eta} $  &  $ n \times r^{2} $, $ n \times r $   &  &   $ \cO(nr^{3}) $  \\
            $  \Sigma_{\odot} = \Sigma_{\circ 2} \otimes \Sigma_{\eta} $  &  $ r^{2} \times r^{2} $, $ r \times r $  & $ \Sigma_{\circ 2} $, $ \Sigma_{\eta} $ sparse diagonal  &  $ \cO(r^{3}) $   \\
            $  V_{\odot} = V_{\circ 2} \kt V_{\eta} $  &  $ n \times r^{2} $, $ n \times r $   &  &  $ \cO(nr^{3}) $  \\
          \bottomrule
      \end{tabular}
   \end{center}
\end{table}

\begin{remark}
We have the following relationships between the discretizations and the derivative/gradient:
\[
   \begin{tikzcd}
      \cF \arrow{r}{\mathrm{discr.}} \arrow[swap]{d}{\mathrm{var.\ der.}}  &  F  \arrow{d}{\mathrm{grad.}} \\
      \frac{\delta \cF}{\delta w}  \arrow{r}{\mathrm{discr.}} &  G
   \end{tikzcd}
\]
i.e., the discretization of $ \cF $ is $ F $, and the Euclidean gradient of $ F $ is $ G $. This is equivalent to computing the variational derivative $ \frac{\delta \cF}{\delta w} $ first, and then discretizing it to obtain $ G $.
\end{remark}

\section{Low-rank formats for the Fisher--KPP equation}\label{app:discretization_FKPP}

\subsection{Cost function}
In low-rank matrix format, we have
\[
   W = U\Sigma V\tr, \qquad W^{(n)} = U^{(n)} \Sigma^{(n)} (V^{(n)})\tr, \qquad  W^{(n-1)} = U^{(n-1)} \Sigma^{(n-1)} (V^{(n-1)})\tr,
\]
where all the $ U $ and $ V $ factors have size $ n $-by-$ r $, while the $ \Sigma $ factors are stored as sparse diagonal $ r $-by-$ r $ matrices. As in~\ref{app:discr_objective} for the Allen--Cahn equation, the square Hadamard power of $ W^{(n)} $ is factorized as $ \big( W^{(n)} \big)^{\circ 2} = U_{\circ 2}\Sigma_{\circ 2}V_{\circ 2}\tr$, where $ U_{\circ 2} $, $ V_{\circ 2} \in \R^{n \times r^{2}} $, and $ \Sigma_{\circ 2} $ is a sparse diagonal $ r^{2} $-by-$ r^{2} $ matrix.

We call the operations $ G_{W} = U\Sigma $, $ G^{(n)} = U^{(n)} \Sigma^{(n)} $, $ G^{(n-1)} = U^{(n-1)} \Sigma^{(n-1)} $, and $ G_{\odot} = U_{\odot} \Sigma_{\odot} $ .
We point out that $ M_{\mathrm{m}}\tr M_{\mathrm{m}} $ is a symmetric sparse banded matrix with bandwidth 2. This implies that the number of nonzero elements is $ 2(n-2) + 2(n-1)+ n = 5n - 6 \ll n^{2} $, which allows for efficient matrix-matrix products.

Refer to Table~\ref{tab:complexities_FKPP_cost} for details on the computational costs for evaluating~\eqref{eq:FKPP_cost_function} in low-rank format.

\begin{table}[htbp]
   \caption{Asymptotic complexities for FKPP cost function.}
   \label{tab:complexities_FKPP_cost}
   \begin{center}
      \resizebox{\textwidth}{!}{%
      \begin{tabular}{c|c|c|c}  
         \toprule
            Product  & Factor sizes &  Notes on structure and storage  &  Cost   \\
          \midrule
          \multicolumn{3}{l}{$ \frac{1}{2} \trace\!\left( G_{W}\tr (M_{\mathrm{m}}\tr M_{\mathrm{m}}) G_{W} \right) $} &  $ \cO(nr^{2}) $  \\
          \midrule
            $ G_{W}\tr (M_{\mathrm{m}}\tr M_{\mathrm{m}}) $  &  $ n \times r $,  $ n \times n $  &  $ M_{\mathrm{m}}\tr M_{\mathrm{m}} $ sparse banded, &  $ \cO(nr) $ \\
                  &    &  $ 5n-6 $ nonzero coefficients  &  \\
          \midrule
          \multicolumn{3}{l}{$ -\trace\!\left( (V\tr V^{(n-1)}) \big((G^{(n-1)})\tr (M_{\mathrm{p}}\tr M_{\mathrm{m}}) \, G_{W} \big) \right) $} &  $ \cO(nr^{2} + r^{3}) $ \\
          \midrule
            $ V\tr V^{(n-1)} $  &  $ n \times r $,  $ n \times r $  &    &  $ \cO(nr^{2}) $ \\
            $ (G^{(n-1)})\tr (M_{\mathrm{p}}\tr M_{\mathrm{m}}) \, G_{W} $ &  $ n \times r $,  $ n \times n $  &  $ M_{\mathrm{p}}\tr M_{\mathrm{m}} $ sparse banded, &  $ \cO(nr^{2}) $ \\
            $ (V\tr V^{(n-1)}) \big((G^{(n-1)})\tr (M_{\mathrm{p}}\tr M_{\mathrm{m}}) \, G_{W} $  &  $ r \times r $,  $ r \times r $  &  $ 5n-6 $ nonzero coefficients  &  $ \cO(r^{3}) $ \\
          \midrule
          \multicolumn{3}{l}{$ -2\deltatime\trace\!\left( \big( (G^{(n)})\tr M_{\mathrm{m}} G_{W} \big) (V\tr \! R_{\omega} V^{(n)}) \right) $} &  $ \cO(nr^{2} + r^{3}) $ \\
          \midrule
            $ (G^{(n)})\tr \! M_{\mathrm{m}} G_{W} $  &  $ n \times r $,  $ n \times n $, $ n \times r $  &  $ M_{\mathrm{m}} $ sparse banded, & $ \cO(nr^{2}) $   \\
            $ V\tr \! R_{\omega} $  &  $ n \times r $,  $ n \times n $  &   $ 3n-2 $ nonzero coefficients  & $ \cO(nr) $   \\
            $ (V\tr \! R_{\omega}) V^{(n)} $ &  $ r \times n $,  $ n \times r $  &  $ R_{\omega} $ sparse diagonal, & $ \cO(nr^{2}) $   \\
            $ \big( (G^{(n)})\tr M_{\mathrm{m}} G_{W} \big) (V\tr \! R_{\omega} V^{(n)}) $  &  $ r \times r $  & $ n $ nonzero coefficients & $ \cO(r^{3}) $  \\
          \midrule
          \multicolumn{3}{l}{$ 2\deltatime\trace\!\left( (G_{\odot}\tr M_{\mathrm{m}} G_{W}) (V\tr R_{\omega} V_{\odot}) \right) $} & $ \cO(nr^{3} + r^{5}) $ \\
          \midrule
            $ G_{\odot}\tr M_{\mathrm{m}} G_{W} $ & $ n \times r^{2} $,  $ n \times n $, $ n \times r $  &  $ M_{\mathrm{m}} $ sparse banded,  & $ \cO(nr^{2}) $   \\
            $ (V\tr \! R_{\omega}) V_{\odot} $  &  $ r \times n $, $ n \times r^{2} $ & $ 3n-2 $ nonzero coefficients  &  $ \cO(nr^{3}) $  \\
            $ (G_{\odot}\tr M_{\mathrm{m}} G_{W}) (V\tr R_{\omega} V_{\odot}) $  &   $ r^{2} \times r $,  $ r \times r^{2} $   &  &  $ \cO(r^{5}) $ \\
         \bottomrule
      \end{tabular}}
   \end{center}
\end{table}

\subsection{Gradient}
The Euclidean gradient of the Fisher--KPP PDE cost function $ F(W) $ is
\[
  G = \left(M_{\mathrm{m}} W - M_{\mathrm{p}} W^{(n-1)} + 2 \deltatime \left(\big( W^{(n)} \big)^{\circ 2} -  W^{(n)}  \right) R_{\omega} \right)\tr \! M_{\mathrm{m}}.
\]
In low-rank format we have $ G = U_{G} \Sigma_{G} V_{G}\tr $, whose factors are
\[
   U_{G} = \left[  \big( M_{\mathrm{m}}\tr M_{\mathrm{m}} \big) U \quad \big( M_{\mathrm{m}}\tr M_{\mathrm{p}} \big) U^{(n-1)} \quad M_{\mathrm{m}}\tr U_{\circ 2} \quad M_{\mathrm{m}}\tr U^{(n)} \right],
\]
\[
   \Sigma_{G} = \blkdiag \! \left( \Sigma, \ -\Sigma^{(n-1)}, \ 2\deltatime \, \Sigma_{\circ 2}, \ -2\deltatime \, \Sigma^{(n)} \right),
\]
\[
    V_{G} = \left[ V \quad V^{(n-1)} \quad R_{\omega} V_{\circ 2} \quad R_{\omega} V^{(n)} \right]\tr.
\]
Table~\ref{tab:complexities_FKPP_gradient} summarizes the asymptotic complexities for the FKPP factorized gradient.

\begin{table}[htbp]
   \caption{Asymptotic complexities for the FKPP gradient.}
   \label{tab:complexities_FKPP_gradient}
   \begin{center}
      \begin{tabular}{c|c|c|c}  
         \toprule
            Product  & Factor sizes &  Notes on structure and storage  &  Cost   \\
         \midrule
         \midrule
            $ \big( M_{\mathrm{m}}\tr M_{\mathrm{m}} \big) U $  &  $ n \times n $,  $ n \times r $  &  $ M_{\mathrm{m}}\tr M_{\mathrm{m}} $ and $ M_{\mathrm{m}}\tr M_{\mathrm{p}} $ sparse banded,  &  $ \cO(nr) $ \\
            $ \big( M_{\mathrm{m}}\tr M_{\mathrm{p}} \big) U^{(n-1)} $  &  $ n \times n $,  $ n \times r $  &  $ 5n-6 $ nonzero coefficients  &  $ \cO(nr) $ \\
            $ M_{\mathrm{m}}\tr U_{\circ 2} $  &  $ n \times n $, $ n \times r^{2} $  &  $ M_{\mathrm{m}} $ sparse banded,  &  $ \cO(nr^{2}) $ \\
            $ M_{\mathrm{m}}\tr U^{(n)} $  &  $ n \times n $, $ n \times r $  &  $ 3n-2 $ nonzero coefficients  &  $ \cO(nr) $ \\
         \midrule
            $ R_{\omega} V_{\circ 2} $  &  $ n \times n $,  $ n \times r^{2} $  &  $ R_{\omega} $ sparse diagonal,  &  $ \cO(nr^{2}) $ \\
            $ R_{\omega} V^{(n)} $  &  $ n \times n $,  $ n \times r $  &  $ n $ nonzero coefficients  & $ \cO(nr) $   \\
         \bottomrule
      \end{tabular}
   \end{center}
\end{table}

\subsection{Hessian}
The discretized Euclidean Hessian is
\[
   \Hessian F(W)[\eta] = \big( M_{\mathrm{m}}\tr M_{\mathrm{m}} \big) \eta.
\]
The factored form of the discretized Euclidean Hessian is
\[
   H_{W}[\eta] = \left[ \big( M_{\mathrm{m}}\tr M_{\mathrm{m}} \big) U_{\eta} \right] \, \left[ S_{\eta} \right] \, \left[ V_{\eta} \right]\tr.
\]
where $ \eta = U_{\eta} S_{\eta} V_{\eta}\tr $ is a tangent vector in $ \mathrm{T}_{W}\cMr $, in a SVD-like format. The only operation needed is the product $ \big( M_{\mathrm{m}}\tr M_{\mathrm{m}} \big) U_{\eta} $, whose cost is $ \cO(nr) $.

\section{Derivation of the preconditioner}\label{app:preconditioner}

As mentioned in Sect.~\ref{sec:RTR}, the tCG trust-region subsolver can be preconditioned with the inverse of~\eqref{eq:prec_candidate}. However, inverting the matrix $ H_{X} $ directly would be too computationally expensive, taking $\cO(n^{6})$ in this case. A suitable preconditioner can be used to solve this problem, thereby reducing the number of iterations required by the tCG solver. This appendix provides the derivation of such a preconditioner.

\subsection{Applying the preconditioner}
In practice, applying the preconditioner in $ X \in \cMr $ means solving (without explicitly inverting the matrix) for $ \xi \in \mathrm{T}_{X}\cM $ the system
\begin{equation}\label{eq:system_hessian_xi_eta}
    H_{X} \vecop(\xi) = \vecop (\eta),
\end{equation}
where $\eta \in \mathrm{T}_{X}\cM $ is a known tangent vector. This equation is equivalent to
\[
    \PTXM(h_{x}^{2}(A \xi + \xi A)) = \eta.
\]
Using definition~\eqref{eq:def_proj} of the orthogonal projector onto $ \TXMr $, we obtain
\[
   P_{U} (A \xi + \xi A) P_{V} + P_{U}^{\perp} (A \xi + \xi A) P_{V} + P_{U} (A \xi + \xi A) P_{V}^{\perp} = \eta,
\]
which is equivalent to the system
\begin{align}\label{eq:system}
\begin{cases}
   P_{U} (A \xi + \xi A) P_{V} = P_{U} \eta P_{V}, \\
   P_{U}^{\perp} (A \xi + \xi A) P_{V} = P_{U}^{\perp} \eta P_{V}, \\
   P_{U} (A \xi + \xi A) P_{V}^{\perp} = P_{U} \eta P_{V}^{\perp}.
\end{cases}
\end{align}
The main difference w.r.t.~\protect{\cite{VandereyckenV_2010}} is that here, in general, the tangent vectors are not symmetric.
Using the matrix representations~\eqref{eq:tan_vec_format_small_param} of the tangent vectors $ \xi $ and $ \eta $ at $ X = U\Sigma V\tr $
\begin{equation*}
   \xi = UM_{\xi}V\tr + \Upxi V\tr + U(\Vpxi)\tr, \qquad \eta = UM_{\eta}V\tr + \Upeta V\tr + U(\Vpeta)\tr,
\end{equation*}
with $M_{\xi} \in \R^{r \times r}$, $\Upxi \in \R^{m\times r}$, $\Vpxi \in \R^{n \times r}$ such that $(\Upxi)\tr U = (\Vpxi)\tr V = 0$.
Analogously, for the tangent vector $\eta$ we have the constraints $M_{\eta} \in \R^{r \times r}$, $ \Upeta \in \R^{m\times r}$, $\Vpeta \in \R^{n \times r}$ such that $(\Upeta)\tr U = (\Vpeta)\tr V = 0$.

After some manipulations (see appendix~\ref{app:preconditioner}), system~\eqref{eq:system} can be written as
\begin{align}\label{eq:system2}
\begin{cases}
   U\tr \! AU M_{\xi} + U\tr \! A \Upxi + M_{\xi} V\tr \! AV + (\Vpxi)\tr \! AV = M_{\eta}, \\
   P_{U}^{\perp} AU M_{\xi} + P_{U}^{\perp} A \Upxi + \Upxi V\tr \! A V =  \Upeta, \\
   M_{\xi} V\tr \! A P_{V}^{\perp} + U\tr \! A U (\Vpxi)\tr + (\Vpxi)\tr \! A P_{V}^{\perp} = (\Vpeta)\tr,
\end{cases}
\end{align}
where $ M_{\xi} $, $ \Upxi $, and $ \Vpxi $ are the unknown matrices.

The solution flow of system~\eqref{eq:system2} is as follows. From the second and the third equations of~\eqref{eq:system2}, we get $ \Upxi $ and $ \Vpxi $ depending on $ M_{\xi} $, then we insert the expressions obtained in the first equation to get $ M_{\xi} $.

We introduce orthogonal matrices $ Q $ and $  \tQ $ to diagonalize $ U\tr \! A U $ and $ V\tr \! A V $, respectively,
\begin{equation}\label{eq:orthogonal_transformations}
   D = QU\tr \! AU Q\tr, \quad \tD =  \tQ V\tr \! AV  \tQ\tr,
\end{equation}
and use them to define the following matrices
\[
   \hU = UQ\tr, \quad \hV = V\tQ\tr, \quad \hMxi = Q M_{\xi} \tQ\tr, \quad \hMeta = Q M_{\eta} \tQ\tr, 
\]
\[
   \hUpxi = \Upxi \tQ\tr, \quad \hVpxi = \Vpxi Q\tr.
\]
With these transformations, the first equation in~\eqref{eq:system2} becomes (see appendix~\ref{app:preconditioner})
\begin{equation}\label{eq:1st_eq_sys}
   D \hMxi + \hMxi \tD + \hU\tr \! A \hUpxi + (\hVpxi)\tr \! A \hV = \hMeta.
\end{equation}
By using the same transformations, we can also rewrite the second equation in~\eqref{eq:system2} as
\[
   P_{U}^{\perp} \big( A \hU \hMxi +  A \hUpxi + \hUpxi \tD \big) =  \hUpeta,
\]
with the condition $ \hU\tr \hUpxi = 0 $.
The $i$th column of this equation is\footnote{We use the MATLAB notation $ (\colon,i) $ to denote the $i$th column extraction from a matrix.}
\[
   P_{U}^{\perp} \, ( A + \widetilde{d}_{i} I ) \, \hUpxi(\colon,i) =  \hUpeta(\colon,i) - P_{U}^{\perp} A \hU \hMxi(\colon,i), \qquad \hU\tr \hUpxi(\colon,i) = 0,
\] 
where $\widetilde{d}_{i}$, for $ i = 1, \ldots, n$, are the diagonal entries of $ \tD $.
We rewrite this equation as a \emph{saddle-point system}
\begin{equation}\label{eq:saddle_pt_sys}
   \begin{bmatrix}
      A + \widetilde{d}_{i} I   &   \hU   \\
                       \hU\tr   &     0
   \end{bmatrix}
   \begin{bmatrix}
      \hUpxi(\colon,i) \\
      y
   \end{bmatrix}
   = 
   \begin{bmatrix}
      \hUpeta(\colon,i) - P_{U}^{\perp} A \hU \hMxi(\colon,i) \\
      0
   \end{bmatrix},
\end{equation}
for all $ y \in \R^{r} $. This saddle-point system can be efficiently solved with the techniques described in Sect.~\ref{sec:saddle_pt_sys}.

Let us define 
\[
   T_{i} \coloneqq 
   \begin{bmatrix}
      A + \widetilde{d}_{i} I   &   \hU   \\
                       \hU\tr   &     0
   \end{bmatrix}, \qquad
   b_{1}^{i} \coloneqq \begin{bmatrix}
      \hUpeta(\colon,i) \\
      0
   \end{bmatrix}, \qquad
   b_{2}^{i} \coloneqq \begin{bmatrix}
      - P_{U}^{\perp} A \hU \hMxi(\colon,i) \\
      0
   \end{bmatrix}.
\]
The solution of~\eqref{eq:saddle_pt_sys} is given by
\[
   \hUpxi(\colon,i) = \big(T_{i}^{-1} b_{1}^{i} + T_{i}^{-1} b_{2}^{i} \big)(1:n),
\]
where the notation $ (1:n) $ means that we only keep the first $ n $ entries of the vector. In other terms, we have
\begin{equation}\label{eq:hUpxi}
   \hUpxi(\colon,i) = \cTi^{-1} \big( \hUpeta(\colon, i) \big) - \cTi^{-1} \big( P_{U}^{\perp} A \hU \big) \hMxi(\colon, i).
\end{equation}
Here, $\cTi^{-1}$ denotes solving for $ \hUpxi(\colon,i) $ the $i$th saddle-point system, corresponding to~\eqref{eq:saddle_pt_sys}.

For the third equation in~\eqref{eq:system2}, we proceed analogously as above. After some manipulations, we obtain the saddle-point system
\begin{equation}\label{eq:saddle_pt_sys_tilde}
   \begin{bmatrix}
      A + d_{i} I   &   \hV   \\
           \hV\tr   &     0
   \end{bmatrix}
   \begin{bmatrix}
      \hVpxi(\colon,i) \\
      z
   \end{bmatrix}
   = 
   \begin{bmatrix}
      \hVpeta(\colon,i) - P_{V}^{\perp} A \hV \hMxi\tr(\colon, i) \\
      0
   \end{bmatrix},
\end{equation}
for all $ z \in \R^{r} $.
The solution is 
\begin{equation}\label{eq:hVpxi}
   \hVpxi(\colon,i) = \cTtildei^{-1} \big( \hVpeta(\colon,i) \big) - \cTtildei^{-1} \big( P_{V}^{\perp} A \hV \big) \hMxi\tr(\colon, i).
\end{equation}
Here, $\cTtildei^{-1}$ denotes solving for $ \hVpxi(\colon,i) $ the $i$th saddle-point system, corresponding to~\eqref{eq:saddle_pt_sys_tilde}.

We now go back to the first equation in~\eqref{eq:system2}, in its form given in~\eqref{eq:1st_eq_sys}.
To treat the term $ \hU\tr \! A \hUpxi $ appearing in~\eqref{eq:1st_eq_sys}, let us define the vectors
\[
   v_{i} = \hU\tr \! A \cTi^{-1} \big( \hUpeta(\colon, i) \big), \qquad   w_{i} = \hU\tr \! A \cTi^{-1} \big( P_{U}^{\perp} A \hU \big) \hMxi(\colon, i).
\]
We emphasize that the vector $ v_{i} $ is known, while $ w_{i} $ is not because $ \hMxi $ is unknown.
With these definitions, and~\eqref{eq:hUpxi}, one can easily verify that
\[
   \hU\tr \! A \hUpxi = \left[ v_{1} - w_{1}, \ldots, v_{r} - w_{r} \right].
\]
Similarly, for treating the term $ \hV\tr \! A \hVpxi $, we define the vectors
\[
   \tv_{i} = \hV\tr \! A \cTtildei^{-1} \big( \hVpeta(\colon, i) \big), \qquad   \tw_{i} = \hV\tr \! A \cTtildei^{-1} \big( P_{V}^{\perp} A \hV \big) \hMxi\tr(\colon, i),
\]
then 
\[
   \hV\tr \! A \hVpxi = \left[ \tv_{1} - \tw_{1}, \ldots, \tv_{r} - \tw_{r} \right].
\]
We insert these two expressions in~\eqref{eq:1st_eq_sys} 
\[
   D \hMxi + \hMxi \tD + \hU\tr \! A \hUpxi + (\hVpxi)\tr \! A \hV = \hMeta,
\]
\[
   D \hMxi + \hMxi \tD + \left[ v_{1} - w_{1}, \ldots, v_{r} - w_{r} \right] + \left[ \tv_{1} - \tw_{1}, \ldots, \tv_{r} - \tw_{r} \right]\tr = \hMeta.
\]
The vectors $ w_{i} $ and $ \tw_{i} $ contain the unknown matrix $ \hMxi $, so we leave them on the left-hand side, while since $ v_{i} $ and $ \tv_{i} $ are known, we move them to the right-hand side.
By letting $ K_{i} \coloneqq \tdi I_{r} - \hU\tr \! A \cTi^{-1} \big( P_{U}^{\perp} A \hU \big) $ and $ \tKi \coloneqq d_{i} I_{r} - \hV\tr \! A \cTtildei^{-1} \big( P_{V}^{\perp} A \hV \big) $ for $ i = 1, \ldots, r $, we get
\begin{equation}\label{eq:system_in_hMxi}
   \left[ K_{1}\hMxi(\colon, 1), \ldots, K_{r}\hMxi(\colon, r) \right] +
   \begin{bmatrix}
      \hMxi(1, :)\,\tK\tr_{1} \\
      \vdots \\
      \hMxi(r, :)\,\tK\tr_{r} 
   \end{bmatrix}
   = R,
\end{equation}
where the matrix on the right-hand side is defined by
\[
   R \coloneqq \hMeta - \left[ v_{1}, \ldots, v_{r} \right] - \left[ \tv_{1}, \ldots, \tv_{r} \right]\tr.
\]
Now, we need to isolate $ \hMxi $ in~\eqref{eq:system_in_hMxi}.
We vectorize the first term on the left-hand side of \eqref{eq:system_in_hMxi} and get
\[
   \vecop\!\left[K_{1}\hMxi(\colon, 1) \quad \cdots \quad K_{r}\hMxi(\colon, r) \right] = \cK \vecop(\hMxi),
\]
where $ \cK = \blkdiag(K_{1}, \ldots, K_{r}) $, a block-diagonal matrix with the $ K_{i} $ on the main diagonal.
Vectorizing the second term on the left-hand side of \eqref{eq:system_in_hMxi}, we obtain
\[
   \vecop\!\left(
   \begin{bmatrix}
      \hMxi(1, :)\,\tK\tr_{1} \\
      \vdots \\
      \hMxi(r, :)\,\tK\tr_{r} 
   \end{bmatrix}\right) = \Pi \vecop\!\left[\tK_{1}\hMxi\tr(1,:) \quad \cdots \quad \tK_{r}\hMxi\tr(r,:) \right] = \Pi \cKtilde \vecop(\hMxi\tr),
\]
where $ \cKtilde = \blkdiag(\tK_{1}, \ldots, \tK_{r}) \in \R^{r^{2} \times r^{2}} $ is a block-diagonal matrix, and $\Pi$ is the \emph{perfect shuffle matrix} defined by $ \vecop(X\tr) = \Pi \vecop(X) $ \protect{\cite{Henderson:1981,VanLoan:2000}}.
Wrapping it up, from \eqref{eq:system_in_hMxi} we obtain the vectorized equation
\begin{equation}\label{eq:small_sys}
   \left[ \cK + \Pi \cKtilde \Pi \right] \vecop(\hMxi) = \vecop(R).
\end{equation}
The matrix $ \cK + \Pi \cKtilde \Pi $ is of size $r^{2}$-by-$r^{2}$, so it can be efficiently inverted if the rank $ r $ is really low.
We solve this equation for $ \hMxi $, and then we use~\eqref{eq:hUpxi} and~\eqref{eq:hVpxi} to find $ \hUpxi $ and $ \hVpxi $, respectively. Finally, undoing the transformations done by $ Q $ and $ \tQ $, we find the components of $\xi$
\[
   M_{\xi} = Q\tr \hMxi \tQ, \qquad \Upxi = \hUpxi \tQ, \qquad \Vpxi = \hVpxi Q,
\]
and thus the tangent vector $ \xi $ such that~\eqref{eq:system_hessian_xi_eta} is satisfied.

\subsection{Efficient solution of the saddle-point system}\label{sec:saddle_pt_sys}
We use a Schur complement idea to efficiently solve the saddle-point problem~\eqref{eq:saddle_pt_sys} and invert the $ T_{i} $. See the techniques described in~\protect{\cite{Benzi:2005}}.
Let
\[
   B_{i} = \left[\hUpeta(\colon,i) \qquad P_{U}^{\perp} A \hU \right] \in \R^{n \times r + 1}.
\]
The system $ \cTi(X) = B_{i} $ can be solved for $X$ by eliminating the (negative) Schur complement
$ S_{i} = \hU\tr (A + \tdi I)^{-1} \hU $. This gives
\begin{equation}\label{eq:Ni}
   N_{i} = S_{i}^{-1}\big(\hU\tr (A + \tdi I)^{-1} B_{i}\big),
\end{equation}
\begin{equation}\label{eq:Xi}
   X_{i} = (A + \tdi I)^{-1} B_{i} - (A + \tdi I)^{-1} \hU N_{i}.
\end{equation}
We use Cholesky factorization to solve the system~\eqref{eq:Ni}. For solving~\eqref{eq:Xi}, we use a sparse solver for $ (A + \tdi I)^{-1} B_{i} $ and $ (A + \tdi I)^{-1} \hU $, for example, MATLAB backslash $ (A + \tdi I) \backslash B_{i} $ and $ (A + \tdi I) \backslash \hU $. Equation~\eqref{eq:small_sys} is a linear system of size $ r^{2} $. 

\subsection{Algebraic manipulations}
\subsubsection{First equation}\label{sec:app_eq1}
The first equation in system~\eqref{eq:system} becomes
\begin{equation*}
\begin{split}
   P_{U} (A (UM_{\xi}V\tr & + \Upxi V\tr + U(\Vpxi)\tr) + (UM_{\xi}V\tr + \Upxi V\tr + U(\Vpxi)\tr) A) P_{V} \\
   & = P_{U} \big( UM_{\eta}V\tr + \Upeta V\tr + U(\Vpeta)\tr \big) P_{V}.
\end{split}
\end{equation*}
Using the gauge conditions
\begin{footnotesize}
\begin{equation*}
\begin{split}
   UU\tr \! AU M_{\xi}V\tr & + UU\tr \! A \Upxi V\tr + \cancel{UU\tr \! AU(\Vpxi)\tr VV\tr} + UM_{\xi}V\tr \! A VV\tr + \\
   & + \cancel{UU\tr \Upxi V\tr \! A VV\tr} + U(\Vpxi)\tr \! A VV\tr = UM_{\eta}V\tr + \cancel{UU\tr \Upeta V\tr} + \cancel{U(\Vpeta)\tr VV\tr},
\end{split}
\end{equation*}
\end{footnotesize}
Left-multiplying by $ U\tr $ and right-multiplying by $ V $ we get
\[
   U\tr \! AU M_{\xi} + U\tr \! A \Upxi + M_{\xi} V\tr \! AV + (\Vpxi)\tr \! AV = M_{\eta}.
\]
With the transformations introduced in~\eqref{eq:orthogonal_transformations}, it becomes
\[
   QU\tr \! AU Q\tr Q M_{\xi} \tQ\tr + QU\tr \! A \Upxi \tQ\tr + Q M_{\xi} \tQ\tr\tQ V\tr \! AV \tQ\tr + Q (\Vpxi)\tr \! AV \tQ\tr = Q M_{\eta} \tQ\tr,
\]
\[
   D \hMxi + \hMxi \tD + \hU\tr \! A \hUpxi + (\hVpxi)\tr \! A \hV = \hMeta.
\]
which is the first equation in system~\eqref{eq:system2}.

\subsubsection{Second equation}\label{sec:app_eq2}
Analogously for the second equation in~\eqref{eq:system}, i.e.,
\[
   P_{U}^{\perp} (A \xi + \xi A) P_{V} = P_{U}^{\perp} \eta P_{V},
\]
\begin{equation*}
\begin{split}
   P_{U}^{\perp} (A (UM_{\xi}V\tr + \Upxi V\tr + U(\Vpxi)\tr) & + (UM_{\xi}V\tr + \Upxi V\tr + U(\Vpxi)\tr) A) P_{V} \\
   & = P_{U}^{\perp} \big( UM_{\eta}V\tr + \Upeta V\tr + U(\Vpeta)\tr \big) P_{V}
\end{split}
\end{equation*}
we obtain
\[
   P_{U}^{\perp} AU M_{\xi}V\tr + P_{U}^{\perp} A \Upxi V\tr + \cancel{P_{U}^{\perp} A U(\Vpxi)\tr VV\tr} + P_{U}^{\perp} \Upxi V\tr \! A VV\tr = P_{U}^{\perp} \Upeta V\tr,
\]
then, using $ P_{U}^{\perp} = I - UU\tr $,
\[
   P_{U}^{\perp} AU M_{\xi}V\tr + P_{U}^{\perp} A \Upxi V\tr + \Upxi V\tr \! A VV\tr - \cancel{UU\tr \Upxi V\tr \! A VV\tr} =  \Upeta V\tr - \cancel{UU\tr \Upeta V\tr}.
\]
Right-multiplying by $V$ we get the second equation in system~\eqref{eq:system2}, i.e.,
\[
   P_{U}^{\perp} AU M_{\xi} + P_{U}^{\perp} A \Upxi + \Upxi V\tr \! A V =  \Upeta.
\]
With the transformations introduced in~\eqref{eq:orthogonal_transformations}, this becomes
\[
   P_{U}^{\perp} AU Q\tr Q M_{\xi}\tQ\tr + P_{U}^{\perp} A \Upxi\tQ\tr + \Upxi \tQ\tr\tQ V\tr \! A V\tQ\tr =  \Upeta\tQ\tr,
\]
\[
   P_{U}^{\perp} A \hU \hMxi + P_{U}^{\perp} A \hUpxi + \hUpxi \tD =  \hUpeta,
\]
\[
   P_{U}^{\perp} \big( A \hU \hMxi +  A \hUpxi + \hUpxi \tD \big) =  \hUpeta.
\]

\subsubsection{Third equation}\label{sec:app_eq3}
The third equation in~\eqref{eq:system}
\[
   P_{U} (A \xi + \xi A) P_{V}^{\perp} = P_{U} \eta P_{V}^{\perp},
\]
becomes
\begin{equation*}
\begin{split}
   P_{U} (A (UM_{\xi}V\tr & + \Upxi V\tr + U(\Vpxi)\tr) + (UM_{\xi}V\tr + \Upxi V\tr + U(\Vpxi)\tr) A) P_{V}^{\perp} \\
   & = P_{U} \big( UM_{\eta}V\tr + \Upeta V\tr + U(\Vpeta)\tr \big) P_{V}^{\perp}.
\end{split}
\end{equation*}
\begin{equation*}
\begin{split}
   P_{U} A (UM_{\xi}V\tr & + \Upxi V\tr + U(\Vpxi)\tr) P_{V}^{\perp} +  P_{U} (UM_{\xi}V\tr + \Upxi V\tr + U(\Vpxi)\tr) A P_{V}^{\perp} \\
   & = P_{U} \big( UM_{\eta}V\tr + \Upeta V\tr + U(\Vpeta)\tr \big) P_{V}^{\perp}.
\end{split}
\end{equation*}
Using the gauge conditions
\begin{footnotesize}
\[
   UU\tr \! A U (\Vpxi)\tr -\cancel{UU\tr \! AU (\Vpxi)\tr VV\tr} + UM_{\xi}V\tr \! A P_{V}^{\perp} + \cancel{UU\tr \Upxi V\tr \! A P_{V}^{\perp}} + U (\Vpxi)\tr \! A P_{V}^{\perp} = U (\Vpeta)\tr P_{V}^{\perp}
\]
\end{footnotesize}
Left-multiplying by $ U\tr $ and noting that $ (\Vpeta)\tr P_{V}^{\perp} = (\Vpeta)\tr $, we obtain the third equation in system~\eqref{eq:system2}.
\[
   M_{\xi} V\tr \! A P_{V}^{\perp} + U\tr \! A U (\Vpxi)\tr + (\Vpxi)\tr \! A P_{V}^{\perp} = (\Vpeta)\tr.
\]

\bibliographystyle{alpha_init}

\begin{footnotesize}
   \bibliography{PrecRTR_biblio.bib}
\end{footnotesize}

\end{document}